\documentclass[10pt,a4paper,oneside,final]{amsart}
\usepackage[T1]{fontenc}
\usepackage[utf8]{inputenc}
\usepackage{lmodern}
\usepackage[final]{microtype}

\usepackage[backend=biber,giveninits=true]{biblatex}

\usepackage{enumitem}
\usepackage{hyperref}

\usepackage{xurl}
\hypersetup{breaklinks=true}

\usepackage{amssymb}
\usepackage{mathtools}
\usepackage{tikz-cd}

\usepackage{amsthm}
\theoremstyle{plain}
\newtheorem{theorem}{Theorem}[section]
\newtheorem{proposition}[theorem]{Proposition}
\newtheorem{lemma}[theorem]{Lemma}
\newtheorem*{convention}{Convention}
\newtheorem*{assumption}{Assumption}
\newtheorem{definition}[theorem]{Definition}
\theoremstyle{remark}
\newtheorem{remark}[theorem]{Remark}
\newtheorem{example}[theorem]{Example}

\newtheorem*{sublemma*}{Sublemma}

\newcommand{\wh}[1]{\hat{#1}}

\newcommand{\wt}[1]{\tilde{#1}}
\newcommand{\dprime}{{\prime\prime}}

\newcommand{\SU}{\mathrm{SU}}
\newcommand{\CP}{\mathbb{C}P}
\newcommand{\N}{\mathbb{N}}
\newcommand{\Z}{\mathbb{Z}}
\newcommand{\K}{\mathbb{K}}
\newcommand{\R}{\mathbb{R}}
\renewcommand{\SS}{\mathcal{S}}
\newcommand{\Pp}{\mathrm{P}}
\newcommand{\QQ}{\mathcal{Q}}
\newcommand{\HH}{\mathcal{H}}
\newcommand{\Hh}{\mathrm{H}}
\newcommand{\rmdeg}{\mathrm{deg}}
\newcommand{\Trees}{\Upsilon}
\newcommand{\Filtr}{\mathcal{F}}

\newcommand{\ALambda}{\Theta}

\newcommand{\Om}{\Omega}
\newcommand{\HDR}{\mathrm{H}_{\mathrm{dR}}}
\newcommand{\Dd}{\mathrm{d}}

\newcommand{\Or}{\varepsilon}
\newcommand{\la}{\langle}
\newcommand{\ra}{\rangle}

\DeclareMathOperator{\Ev}{ev}
\DeclareMathOperator{\im}{im}
\DeclareMathOperator{\Span}{span}
\newcommand{\Id}{{{\mathchoice {\rm 1\mskip-4mu l} {\rm 1\mskip-4mu l} {\rm 1\mskip-4.5mu l} {\rm 1\mskip-5mu l}}}}
\DeclareMathOperator{\Char}{\mathrm{char}}
\newcommand{\onto}{\twoheadrightarrow}
\newcommand{\into}{\hookrightarrow}

\begin{document}

\title{Hodge decompositions and differential Poincar\'e duality models}
\author{Pavel H\'ajek}

\begin{abstract}
We extend a CDGA $V$ with a perfect pairing of degree~$n$ on cohomology to a CDGA $\wh V$ with a pairing of degree $n$ on chain level such that~$\wh V$ admits a Hodge decomposition and retracts onto $V$ preserving the pairing on cohomology; here we suppose that $V$ is either 1-connected, or that $V$ is connected, of finite type, and $n$ is odd. 
We show that a Hodge decomposition of~$\wh V$ induces a differential Poincar\'e duality model of~$V$ in a natural way.
Assuming that $\Hh(V)$ is 1-connected, we apply our extension to a Sullivan model of $V$ in the proof of the existence and ``uniqueness'' of a 1-connected differential Poincar\'e duality model of~$V$ by Lambrechts \& Stanley; we eliminate their extra assumptions in the uniqueness statement, including $\Hh^2(V)=0$ if $n$ is odd.
\end{abstract}
\maketitle
\tableofcontents
\section{Introduction}\label{Sec:1}

A~\emph{differential Poincar\'e duality algebra} (\emph{dPD algebra}) $(V,\Dd,\Or)$ of degree $n\in\N_0$ over a field~$\K$ is a non-negatively graded unital commutative differential graded algebra (\emph{CDGA}) $(V=\oplus_{i\in\N_0} V^i,\Dd)$ over~$\K$ together with an \emph{orientation} $\Or\colon V\to \K$ of degree~$n$---that is, a linear map $\Or\colon V\to \K$ such that $\im\Dd\subset\ker\Or$, $\Or|_{\ker\Dd}\neq 0$, and $V^i\subset\ker\Or$ whenever $i\neq n$---such that the pairing $\la-,-\ra\colon V\times V\to \K$ given by $\la v_1,v_2\ra \coloneqq \Or(v_1 v_2)$ is perfect.
We also say that $(V,\Dd,\Or)$ satisfies \emph{Poincar\'e duality}.

A~\emph{Poincar\'e DGA} (\emph{PDGA}) $(V,\Dd,\Or_*)$ of degree $n$ over $\K$ is a non-negatively graded unital CDGA $(V,\Dd)$ over $\K$ whose \emph{cohomology} $\Hh(V)\coloneqq \Hh(V,\Dd)$ is equipped with an orientation $\Or_*\colon\Hh(V) \rightarrow \K$ of degree $n$ that makes it into a \emph{Poincar\'e duality algebra} (\emph{PD algebra})---that is, a dPD algebra with $\Dd\equiv 0$.
An \emph{oriented PDGA} $(V,\Dd,\Or)$ of degree $n$ over~$\K$ is a PDGA $(V,\Dd,\Or_*)$ of degree $n$ over~$\K$ together with an orientation $\Or\colon V\to\K$ of degree $n$ on $(V,\Dd)$ that induces the orientation $\Or_*\colon\Hh(V)\to\K$. 
We refer to $\Or$ as \emph{chain-level orientation} and to $\Or_*$ as \emph{cohomology orientation}.
Notice that a PDGA is only required to satisfy Poincar\'e duality on cohomology, and that a dPD algebra is an oriented PDGA such that Poincar\'e duality holds also on chain level.
An archetypal example of an oriented PDGA of degree $n$ that is \emph{not} a dPD algebra is the \emph{oriented de Rham algebra} $(\Om(X),\Dd,\Or_X)$ of an oriented closed $n$-manifold $X$ with the orientation given by the integration $\Or_X(\omega)\coloneqq\int_X \omega$ for $\omega\in\Om(X)$.

A morphism $f\colon V\to V^\prime$ of PDGAs $(V,\Dd,\Or_*)$ and $(V^\prime,\Dd^\prime,\Or^\prime_*)$ of the same degree~$n$ (\emph{PDGA morphism}) is a morphism of the underlying DGAs (DGA morphism) that preserves cohomology orientation.
We say that $(V,\Dd,\Or_*)$ and $(V^\prime,\Dd^\prime,\Or^\prime_*)$ are \emph{weakly homotopy equivalent as PDGAs} if they can be connected by a zig-zag of PDGA quasi-isomorphisms.
Any dPD algebra $(V^\prime,\Dd^\prime,\Or^\prime)$ that is weakly homotopy equivalent to $(V,\Dd,\Or_*)$ as a PDGA is called a \emph{differential Poincar\'e duality model} (\emph{dPD model}) of $(V,\Dd,\Or_*)$.
The following fundamental results were proven in~\cite{Lambrechts2007}: 

\begin{enumerate}
	\item[(E)] \emph{Existence statement} \cite[Theorem~1.1]{Lambrechts2007}: If $\Hh(V)$ is 1-connected, then there is a dPD algebra $(V^\prime,\Dd^\prime,\Or^\prime)$ that is weakly homotopy equivalent to $(V,\Dd,\Or_*)$ as a DGA---that is, $(V,\Dd)$ and $(V^\prime,\Dd^\prime)$ are connected by a zig-zag of DGA quasi-isomorphisms.
	Such $(V^\prime,\Dd^\prime,\Or^\prime)$ is called a dPD model in \cite{Lambrechts2007}.
	\item[(U)] \emph{Uniqueness statement} \cite[Theorem~7.1]{Lambrechts2007}: Suppose that~$(V,\Dd,\Or)$ and~$(V^\prime,\Dd^\prime,\Or^\prime)$ are dPD algebras of degree $n$ that are weakly homotopy equivalent as DGAs.
	Suppose additionally that $n\ge 7$, that $V$ and $V^\prime$ are 2-connected, and that $\Hh(V)\simeq \Hh(V^\prime)$ is $3$-connected. Then there is a dPD algebra $(V^\dprime,\Dd^\dprime,\Or^\dprime)$ of degree~$n$ and DGA quasi-isomorphisms $f\colon V\to V^\dprime$ and $f^\prime\colon V^\prime\to V^\dprime$.
\end{enumerate}
Given an oriented PDGA $(V,\Dd,\Or)$, we define the \emph{degenerate subspace} $V_\rmdeg\coloneqq\{v\in V\mid v\perp V\}$ and the \emph{nondegenerate quotient} $\QQ(V)\coloneqq V/V_\rmdeg$, where $\perp$ denotes the perpendicularity relation with respect to $\la-,-\ra$.
The nondegenerate quotient $\QQ(V)$ is naturally a CDGA such that the quotient map $\pi_\QQ\colon V\onto \QQ(V)$ is a DGA morphism, and there is a unique orientation $\Or_\QQ\colon \QQ(V)\to\K$ such that $\Or_\QQ \circ \pi_\QQ = \Or$.
One can easily see that~$\pi_\QQ\colon V\onto \QQ(V)$ is a quasi-isomorphism if and only if $\ker\pi_\QQ = V_\rmdeg$ is acyclic, and that $\QQ(V)$ is a dPD algebra provided that $V$ is of finite type.

The pivotal piece of \cite{Lambrechts2007} is a construction of an oriented PDGA $(\wh V,\wh \Dd,\wh \Or)\supset (V,\Dd,\Or)$ such that that both the inclusion $\iota\colon V\into \wh V$ and the quotient map $\pi_\QQ\colon \wh V \onto \QQ(\wh V)$ are DGA quasi-isomorphisms.
Their construction proceeds iteratively in degrees $k=\lceil n/2\rceil, \dotsc,n$ by adjoining new algebraic generators that ``kill'' some chosen generators of $\Hh^k(V_\rmdeg)$---the so called ``orphans''---but leave $\Hh^k(V)$ unchanged.
This method works if $V$ is 1-connected; hence, in~\cite{Lambrechts2007}, they apply it to a minimal Sullivan model $\rho\colon \Lambda U \to V$, which is 1-connected and of finite type provided that $\Hh(V)$ is.
The nondegenerate quotient $\QQ(\wh{\Lambda U})$ is then a 1-connected dPD model of $V$. 

They observed in~\cite{Van2019} that if $V$ admits a \emph{Hodge decomposition}---that is, a direct sum decomposition $V=\HH\oplus \im\Dd\oplus C$, where $\HH$ is a complement of $\im\Dd$ in $\ker\Dd$ and~$C$ a complement of $\ker\Dd$ in~$V$ such that $C\perp \HH\oplus C$---then $\Hh(V_\rmdeg)=0$.
They call an oriented PDGA $(V,\Dd,\Or)$ that admits a Hodge decomposition \emph{of Hodge type}.
By a result of~\cite{Cieliebak2015}, every \emph{cyclic cochain complex}---that is, a cochain complex with a perfect pairing compatible with $\Dd$---is of Hodge type, provided that $\Char(\K)\neq 2$.
This can be used to prove the reverse implication---if $\Hh(V_\rmdeg)=0$, then a Hodge decomposition exists---provided that $V$ is of finite type.
Therefore, in this case, $\pi_\QQ\colon V\onto \QQ(V)$ is a quasi-isomorphism if and only if~$V$ is of Hodge type.

Given this knowledge, we propose an alternative construction of $(\wh V,\wh \Dd, \wh \Or)$ that focuses on removing the obstruction to the existence of a Hodge decomposition rather than on ``killing the orphans''; we will refer to our $\wh V$ as \emph{extension of Hodge type}.
We start with a \emph{pre-Hodge decomposition} $V = \HH \oplus \im\Dd\oplus C$, which is defined as a Hodge decomposition without the assumption $C\perp C$, and which exists for any PDGA.
In order to achieve $C\perp C$, equivalently $\QQ(C)=0$, we proceed iteratively in degrees $k=\lceil n/2\rceil, \dotsc, n$ as follows.
We choose linearly independent elements $q_\alpha\in C^k$ ($\alpha\in I$) whose images under $\pi_\QQ$ generate the vector space $\QQ^{k}(C)$, we introduce new algebraic generators $x_\alpha$, $\Dd x_\alpha$ with $\deg \Dd x_\alpha = \deg x_\alpha + 1 = k$ for each $\alpha\in I$, and we extend the pairing such that $\la q_\alpha, c \ra = - \la \Dd x_\alpha, c \ra$ for all $c\in C^{n-k}$.
Replacing $q_\alpha\mapsto q_\alpha + \Dd x_\alpha$ in~$C^k$ for all $\alpha\in I$, we obtain a complement~$\wh C^k$ of $(\ker\wh\Dd)^k$ in~$\wh V^k$ such that $\wh C^k\perp \wh{C}^{n-k}\oplus \HH^{n-k}$.
This method works if~$V$ is 1-connected.
The difficulty for $V^1\neq 0$ is that $\Span_\K\{ \gamma x_\alpha \mid \gamma\in V^1, \alpha\in I\}\subset \wh V^{k}$ could possibly induce elements in $\QQ^{k}(\wh C)$ that do not come from $\QQ^{k}(C)$, so we may have to adjoin another set of generators in order to deal with them. 
Nonetheless, we observed that this cannot happen indefinitely if $V$ is of finite type and $n\neq 2k$, which is due to Poincar\'e duality $\QQ^k(\wh C)\simeq \QQ^{n-k}(\wh C)$ and the fact that $\dim \QQ^{n-k}(\wh C)\le \dim C^{n-k}$ for $n-k<k$.
This allows us to extend our construction to the connected case. 

We proposed our construction in~\cite{MyPhD} after reading \cite{Van2019} and committing to obtain the results (E) and (U) of~\cite{Lambrechts2007} in the formalism of Hodge decompositions.
We have not studied the relation of our construction to the construction in \cite{Lambrechts2007} in detail.
However, one obvious difference is that our differential on $\wh V$ is the tensor product differential (which allows us to talk about retractions later), whereas their differential is twisted.

The following theorems correspond to Propositions~\ref{Prop:ExtensionSimplyConnected} and~\ref{Prop:ExtensionConnected}, respectively.
The proofs comprise a detailed analysis of our extension of Hodge type.

\begin{theorem}[1-Connected case]\label{Thm:1}
	Let $(V,\Dd,\Or)$ be a \emph{1-connected} oriented PDGA of degree $n\in\N_0$ over a field $\K$.
If $\Char(\K)=2$ and $n\ge 4$ is even, we additionally assume that $V^{n/2}$ admits a Hodge decomposition.
Then there exists a connected oriented PDGA $(\wh V,\wh \Dd,\wh \Or)\supset (V,\Dd,\Or)$ of degree $n$ of Hodge type and a \emph{PDGA retraction} $\pi\colon \wh V\onto V$---that is, a PDGA quasi-isomorphism $\pi\colon \wh V\onto V$ such that $\pi|_V = \Id_V$.
Moreover, $\wh V$ is $1$-connected for $n\neq 4$, and $\wh V$ is of finite type provided that~$V$ is.
\end{theorem}

\begin{theorem}[Connected case of finite type]\label{Thm:2}
Let $(V,\Dd,\Or)$ be a \emph{connected} oriented PDGA of degree $n\in\N_0$ \emph{of finite type} over a field $\K$.
If $n$ is even, we additionally assume that~$V^{n/2}$ admits a Hodge decomposition.
Then there exists a connected oriented PDGA $(\wh V,\wh \Dd,\wh \Or)\supset (V,\Dd,\Or)$ of degree $n$ of Hodge type of finite type, and there is a PDGA retraction $\pi\colon \wh V\onto V$.
\end{theorem}

The necessity of the additional assumption of Theorem~\ref{Thm:1} for $\Char(\K)=2$ is demonstrated in Example~\ref{Ex:NoHodgeExt}.
However, this does not contradict~\cite{Lambrechts2007} because a PDGA whose degenerate subspace is acyclic is not necessarily of Hodge type if $\Char(\K)=2$.
The necessity of the additional assumption $n\neq 4$ for 1-connectedness is demonstrated in Example~\ref{Ex:Degree4}.

The following theorems are our versions of (E) and (U) and correspond to Propositions~\ref{Prop:Existence} and \ref{Prop:Uniqueness}, respectively.
Our proofs are almost identical to the proofs in~\cite{Lambrechts2007}, except that we use our extension of Hodge type.
However, unlike~\cite{Lambrechts2007}, we prefer to assume that $\Char(\K)=0$ because our reference on Sullivan algebras does so.
Nevertheless, we discuss in Remark~\ref{rem:assumptiononchar} that only the assumption $\Char(\K)\neq 2$ for even $n$ due to Theorem~\ref{Thm:1} should be necessary for our approach via Hodge decompositions to work.

\begin{theorem}[Existence]\label{Thm:3}
	Every PDGA $(V,\Dd,\Or_*)$ of degree $n\in\N_0$ with 1-connected cohomology admits a 1-connected dPD model $(V^\prime,\Dd^\prime,\Or^\prime)$ of the form 
	\begin{equation*}
		 \begin{tikzcd}
			 & \wh{\Lambda U}\arrow["\rho\,\circ\,\pi"']{dl}\arrow[two heads,"\pi_\QQ"]{dr} & \\
			 V & & V^\prime \coloneqq \QQ(\wh{\Lambda U}),
		 \end{tikzcd}
	 \end{equation*}
	 where $\rho\colon \Lambda U\to V$ is a minimal Sullivan model of $(V,\Dd)$, $\wh{\Lambda U}\supset \Lambda U$ is a 1-connected oriented PDGA of Hodge type of finite type, $\pi\colon \wh{\Lambda U}\onto\Lambda U$ is a PDGA retraction, and $\pi_\QQ\colon \wh{\Lambda U}\onto \QQ(\wh{\Lambda U})$ is the quotient map.
\end{theorem}

\begin{theorem}[Uniqueness]\label{Thm:4}
Let $(V,\Dd,\Or)$ and $(V^\prime,\Dd^\prime,\Or^\prime)$ be 1-connected dPD algebras of degree $n\in\N_0$ that are weakly homotopy equivalent as PDGAs.
If $\Hh^2(V) \simeq \Hh^2(V^\prime) =0$, then there is a 1-connected dPD algebra~$(V^\dprime,\Dd^\dprime,\Or^\dprime)$ and PDGA quasi-isomorphisms
\begin{equation*}
\begin{tikzcd}
& V^\dprime & \\
V \arrow[hook]{ur}{\iota}  & & V^\prime \arrow[hook',swap]{ul}{\iota^\prime},
\end{tikzcd}
\end{equation*}
which are injective and preserve chain-level orientation.
\end{theorem}

Theorem~\eqref{Thm:4}, as it is stated, \emph{does not directly generalize} (U): on one hand, we have the stronger assumption that $V$ and $V^\prime$ are weakly homotopy equivalent \emph{as PDGAs}; on the other hand, we have the stronger assertion that $\iota$ and $\iota^\prime$ are injective and preserve chain-level orientation (this holds, in fact, for any PDGA morphism such that the chain-level pairing on the domain is nondegenerate).
A~generalization of (U) can be obtained as follows: we start with a zig-zag of DGA quasi-isomorphisms between~$V$ and~$V^\prime$ and define a new orientation $\wt\Or\colon V\to \K$ on $(V,\Dd)$ such that the isomorphism $\Hh(V)\simeq \Hh(V^\prime)$ induced by the zig-zag intertwines $\wt\Or_*$ and $\Or^\prime_*$.
We have $V^{n}\simeq V^{0}\simeq \K$ by Poincar\'e duality, so $\wt\Or$ must be a nonzero multiple of $\Or^\prime$; hence, $(V,\Dd,\wt\Or)$ is a dPD algebra.
We can now apply Theorem~\ref{Thm:4} to $(V,\Dd,\wt\Or)$ and $(V^\prime,\Dd^\prime,\Or^\prime)$, obtaining DGA quasi-isomorphisms $\iota\colon V\to V^\dprime$ and~$\iota^\prime\colon V^\prime\to V^\dprime$ as in~(U).
Moreover, both are injective, but only $\iota^\prime$ preserves the original orientation.

A conjecture at the end of~\cite{Lambrechts2007} states that the additional assumptions of~(U) can be dropped.
Using our extension of Hodge type in the connected case, we can show the following:

\begin{theorem}\label{Thm:5}
	If $n$ is \emph{odd}, then Theorem~\ref{Thm:4} holds without the assumption $\Hh^2(V) \simeq \Hh^2(V^\prime) =0$ provided that we allow $V^{\dprime 1}\neq 0$.
\end{theorem}

The proof is identical to the proof of Theorem~\ref{Thm:4} and relies on the construction of an extension of Hodge type of the relative Sullivan algebra $V\otimes V^\prime \otimes \Lambda(U[1])$ (a~replacement of the pushout of $V\leftarrow\Lambda U\rightarrow V$).
This is now only connected since $U[1]^1\simeq \Hh^2(V)$, so Theorem~\ref{Thm:2} has to be used. 

Our motivation to investigate Hodge decompositions and dPD models comes from the study of a natural \emph{IBL$_\infty$ structure} on cyclic Hochschild cochains of the \emph{de Rham cohomology} $\HDR(X)\coloneqq \Hh(\Om(X),\Dd)$ of a closed oriented manifold $X$.
This structure was proposed in~\cite{Cieliebak2015} and is related to symplectic geometry of the cotangent bundle~$T^*X$.
It can be obtained by twisting a canonical IBL structure with a natural \emph{Maurer-Cartan element} $\mathfrak{m}$ that formally corresponds to the homotopy transfer of the structure of a DGA with a pairing on $\Om(X)$ to the structure of a \emph{quantum A$_\infty$-algebra} on $\HDR(X)$---that is, a generalization of \emph{cyclic~A$_\infty$ algebra} with operations in each genera.
The Maurer-Cartan element is constructed in~\cite{chernsimons} as a sum of \emph{Feynman integrals} associated to ribbon graphs with trivalent vertices decorated with integration variables $x_i$, univalent vertices (adjacent to $x_i$) decorated with harmonic forms~$h_j(x_i)$, and edges (connecting $x_i$ and $x_j$) decorated with the integral kernel $\Pp(x_i,y_j)$ of a propagator $\Pp\colon \Om(X)\to\Om(X)$.
Notice that this data is similar to the data used in the construction of an effective action in \emph{perturbative Chern-Simons theory} (see~\cite{CattaneoMnev}).
In~\cite{chernsimons} they show that $\mathfrak{m}$ is well-defined and only depends on $\Pp$ up to gauge equivalence, which induces a homotopy of IBL$_\infty$ algebras.

In~\cite{triplepaper} we use Theorems~\ref{Thm:3} and~\ref{Thm:4} to show that if $\HDR(X)$ is 2-connected, then the twisted IBL$_\infty$ structure is homotopy equivalent to a canonical dIBL structure on cyclic Hochschild cochains of a dPD model of $(\Om(X),\Dd,\Or_{X*})$.
This was conjectured in~\cite{MyPhD}.
If $n$ is odd, we can use Theorem~\ref{Thm:5} and immediately extend this result to the case that~$\HDR(X)$ is only 1-connected.
Note that if $X$ is 1-connected, then the canonical dIBL structure on cyclic Hochschild cochains of a dPD model of $(\Om(X),\Dd,\Or_{X*})$ is a chain model of the \emph{equivariant Chas-Sullivan string topology} of~$X$.
This was shown in~\cite{Naef2019} using a dPD model for the configuration space of two points.

Hodge decompositions were also used in~\cite{LazarevHodge} for the construction of \emph{minimal models of algebras over cyclic operads} (in~\cite{LazarevFeynman} they continue this program by constructing minimal models for algebras over Feynman transform of modular operads).

\vspace{1ex}

\textit{Acknowledgments:}
I~thank Prof.\,Dr.\,H{\^{o}}ng V{\^{a}}n L{\^{e}} for suggesting that the existence of a Hodge decomposition and $\Hh(V_\rmdeg)=0$ are equivalent and for her interest in my work.
I~thank Thorsten Hertl for showing me how to compute the cohomology ring of a connected sum.
I~thank Prof.\,Dr.\,Janko Latschev and Dr.\,Evgeny Volkov for reading the text and suggesting improvements.
I~additionally thank Janko for his patience with me postponing finishing the paper.
I thank Inessa Tunieva for motivating me to finally finish the paper and for discussing the Introduction.
I~thank Prof.\,Dr.\,Kai Cieliebak for his supervision during the writing of my PhD~thesis, which this work stems from, and for his continuous support.
I~thank University of Hamburg, University of Augsburg, and Mittag-Leffler Institute for their support and hospitality.
Finally, I thank the anonymous referee for his/her suggestions and patience. 

\section{Elementary notions}\label{Sec:2}

We write $\N\coloneqq\{1,2,\dotsc\}$ for the set of natural numbers, and we set $\N_0\coloneqq \N\cup\{0\}$.
The symbol $\K$ stands for a field.
Given a real number $\lambda\in \R$, we denote by $\lceil \lambda \rceil\in\Z$ the least integer that is greater than or equal to~$\lambda$.

We work with \emph{non-negatively graded vector spaces} $V=\oplus_{i\in\N_0} V^i$ over~$\K$; we set $V^{-i}\coloneqq 0$ for all $i\in\N$.
Morphisms $f\colon V\to V^\prime$ are \emph{graded linear maps}---that is, $f\colon V\to V^\prime$ is $\K$-linear and there is a $d\in\Z$ such that $f(V^i)\subset V^{i+d}$ for all $i\in\Z$; we denote $\deg f \coloneqq d$ and call it the \emph{degree} of~$f$.
All constructions (like direct sum $V_1\oplus V_2$, tensor product $V_1\otimes V_2$, or dual $V^*$) and relations (like being complementary or being a subspace $V_1\subset V_2$) shall be understood in the \emph{graded sense} (degreewise).
A~vector $v\in V$ is called \emph{homogenous} if there is a $d\in\Z$ such that $v\in V^d$; we denote $\deg v\coloneqq d$ and call it the \emph{degree} of~$v$.
A~graded linear map $f\colon V\to V^\prime$ is called \emph{homogenous} if $\deg f=0$.
We say that $V$ is \emph{of finite type} if $\dim V^i<\infty$ for all $i\in\Z$.
Given an interval $I\subset \R$, we denote $V^I\coloneqq \oplus_{i\in I} V^i$, where the sum is over all integers $i\in I$.

We usually refer to a non-negatively graded vector space simply as ``vector space'' or ``graded vector space'' and to a graded linear map simply as ``map'' or ``linear map'' since the meaning is clear from the context.
We usually assume that vectors are homogenous without explicitly stating it.
We write $i\in I$ 

We say that a \emph{bilinear form} $\la-,-\ra \colon V\times V\rightarrow \K$ has \emph{degree} $n\in\N_0$ if $\la v_1,v_2\ra \neq 0$ implies $\deg v_1 + \deg v_2 = n$ for all $v_1, v_2\in V$ (this is equivalent to $\deg{\la-,-\ra}=-n$ for the corresponding linear map $\la-,-\ra\colon V\otimes V\to \K$).
We say that $\la-,-\ra$ is \emph{graded symmetric} if for all $v_1, v_2\in V$ we have $\la v_1, v_2 \ra = (-1)^{\deg v_1 \deg v_2} \la v_2, v_1\ra$.
Suppose that this is the case.
For $v_1,v_2\in V$ we write $v_1\perp v_2$ if $\la v_1,v_2\ra =0$, and for $V_1, V_2\subset V$ we write $V_1\perp V_2$ if $v_1\perp v_2$ holds for all $v_1\in V_1$, $v_2\in V_2$.
We say that $\la-,-\ra$ is \emph{nondegenerate} if $v_1\perp v_2$ for all $v_2\in V$ implies $v_1=0$.
We say that $\la-,-\ra$ is \emph{perfect} if the map $v\in V \mapsto \la v,-\ra\in V^*$ is an isomorphism.
Because~$V$ is non-negatively graded, the nondegeneracy of $\la-,-\ra$ implies $V = V^0\oplus\dotsb\oplus V^n$, so $\la-,-\ra$ is perfect if and only it is nondegenerate and~$\dim V<\infty$.

A \emph{cochain complex} $(V,\Dd)$ is a non-negatively graded vector space $V$ together with a graded linear map $\Dd\colon V\to V$ with $\deg{\Dd}=1$ such that $\Dd\circ\Dd=0$.
We denote by $\Hh(V)\coloneqq \Hh(V,\Dd)$ the \emph{cohomology} of $(V,\Dd)$.
A \emph{differential graded algebra (DGA)} is a cochain complex $(V,\Dd)$ together with an associative product ${\cdot}\colon V\times V\to V$ of degree~$0$ such that the \emph{Leibniz identity} $\Dd (v_1 \cdot v_2) = (\Dd v_1)\cdot v_2 + (-1)^{\deg v_1} v_1\cdot (\Dd v_2)$ holds for all $v_1,v_2\in V$; we usually omit $\cdot$ from the notation.
A~\emph{commutative DGA (CDGA)} must additionally satisfy $v_1\cdot v_2 = (-1)^{\deg v_1 \deg v_2} v_2\cdot v_1$ for all $v_1, v_2\in V$.
A~\emph{unital DGA} contains an element $1\in V^0$, called the \emph{unit}, such that $1\cdot v = v\cdot 1 = v$ for all $v\in V$.
A unital DGA $(V,\Dd)$ is called \emph{connected} if $V^0=\Span_\K\{1\}$ and \emph{$k$-connected} ($k\in \N$) if it is connected and $V^1=V^2=\dotsb=V^k=0$ (another~common term for $1$-connected is \emph{simply connected}).

\section{Pairings and Hodge decompositions}\label{Sec:3}

A \emph{pairing of degree $n\in\N_0$} on a cochain complex $(V,\Dd)$ is a bilinear form $\la-,-\ra\colon V\times V\to \K$ of degree $n$ that for all $v_1, v_2\in V$ satisfies
\begin{equation}\label{Eq:Differential}
	\la \Dd v_1, v_2 \ra = (-1)^{\deg v_1 + 1} \la v_1, \Dd v_2\ra. 
\end{equation}
A pairing on a DGA must for all $v_1, v_2, v_3\in V$ additionally satisfy 
\begin{equation}\label{Eq:Product}
\la v_1 v_2, v_3 \ra = \la v_1, v_2  v_3\ra. 
\end{equation} 
We call the triple $(V,\Dd,\la-,-\ra)$ a \emph{cochain complex (or a DGA) with a pairing of degree~$n$}.
If $\la-,-\ra$ is graded symmetric, then~\eqref{Eq:Differential} is for each $v_1, v_2\in V$ equivalent to
\begin{equation}\label{Eq:DifferentialCyc}
	\la \Dd v_1, v_2 \ra  = (-1)^{1+\deg v_1 \deg v_2} \la \Dd v_2, v_1 \ra,
\end{equation}
and the condition \eqref{Eq:Product} for all $v_1, v_2, v_3\in V$ is equivalent to the condition
\begin{equation} \label{Eq:ProductCyc}
	\la v_1 v_2, v_3 \ra = (-1)^{\deg v_3 (\deg v_1 + \deg v_2)} \la v_3 v_1,v_2\ra.
\end{equation}
If $(V,\Dd)$ is unital, then the graded symmetry of $\la-,-\ra$ follows from~\eqref{Eq:ProductCyc}.
If $(V,\Dd)$ is unital and commutative, then the graded symmetry of $\la-,-\ra$ follows from~\eqref{Eq:Product}.
We see that a bilinear form $\la-,-\ra\colon V\times V\to\K$ on a unital CDGA $(V,\Dd)$ satisfies conditions~\eqref{Eq:Differential} and~\eqref{Eq:Product} if and only if it satisfies conditions~\eqref{Eq:DifferentialCyc} and~\eqref{Eq:ProductCyc}; it is then automatically graded symmetric.

The following notion is used in~\cite{Cieliebak2015}:

\begin{definition}\label{Def:Cyclic}
	A \emph{cyclic cochain complex (or a DGA)} of degree $n\in\N_0$ is a finite-dimensional cochain complex (or a DGA) $(V,\Dd)$ with a nondegenerate graded symmetric pairing $\la-,-\ra\colon V\times V\to \K$ of degree $n$.
\end{definition}

The following notions are used in~\cite{Lambrechts2007}: 

\begin{definition}
An \emph{orientation} of degree $n\in\N_0$ on a cochain complex (or a DGA) $(V,\Dd)$ is a linear map $\Or\colon V \to \K$ that satisfies the following:
\begin{enumerate}
	\item $\Or(v)\neq 0$ implies $\deg v = n$ for all $v\in V$ (equivalently $\deg \Or = -n$),
	\item $\Or \circ \Dd = 0$,
	\item the induced map on cohomology $\Or_*\colon \Hh(V) \to \K$, $[v]\mapsto \Or(v)$ is nonzero.
\end{enumerate}
The triple $(V,\Dd,\Or)$ is called an \emph{oriented cochain complex (or a DGA) of degree $n$}. 
\end{definition}

\begin{definition}\label{Def:dPD}
	A \emph{differential Poincar\'e duality algebra (dPD algebra)} $(V,\Dd,\Or)$ of degree $n\in\N_0$ is an oriented unital CDGA $(V,\Dd,\Or)$ of degree $n$ such that the pairing $\la-,-\ra\colon V\times V\to\K$ that is for all $v_1, v_2\in V$ defined by
	\begin{equation}\label{Eq:pairing-orientation}
		\la v_1, v_2\ra \coloneqq \Or(v_1\cdot v_2)	
	\end{equation}
	is perfect.
	A \emph{Poincar\'e duality algebra (PD algebra)} $(V,\Or)$ of degree $n$ corresponds to a dPD algebra of degree $n$ with $\Dd \equiv 0$.
\end{definition}

Equation~\eqref{Eq:pairing-orientation} on a unital DGA $(V,\Dd)$ defines a one-to-one correspondence of linear maps $\Or\colon V\to \K$ with $\deg\Or=-n$ that satisfy $\Or\circ\Dd=0$ and pairings $\la-,-\ra\colon V\times V\to \K$ of degree~$n$.
An immediate consequence is the following:

\begin{lemma}\label{Lem:CyclicEqualsdPD}
	The notions of a dPD algebra $(V,\Dd,\Or)$ and a unital cyclic CDGA $(V,\Dd,\la-,-\ra)$ agree under the correspondence~\eqref{Eq:pairing-orientation}.\hfill\qed
\end{lemma}

The following notion from~\cite{Van2019} is slightly generalized by allowing $V^i\neq 0$ for $i>n$:

\begin{definition}\label{Def:PDGA}
	A~\emph{Poincar\'e DGA (PDGA)} $(V,\Dd,\Or_*)$ of degree $n\in\N_0$ is a unital CDGA $(V,\Dd)$ together with an orientation $\Or_*\colon\Hh(V)\to\K$ of degree $n$ that makes $(\Hh(V),\Or_*)$ into a Poincar\'e duality algebra.
	An \emph{oriented PDGA} $(V,\Dd,\Or)$ of degree $n$ is a PDGA $(V,\Dd,\Or_*)$ of degree $n$ together with an orientation $\Or\colon V\to \K$ of degree~$n$ on $(V,\Dd)$ such that $\Or_*([v])=\Or(v)$ for all $[v]\in \Hh(V)$.
	We call $\Or_*$ \emph{cohomology orientation} and $\Or$ \emph{chain-level orientation}.  
\end{definition}

 We emphasize that a dPD algebra has a perfect pairing on chain level whereas a PDGA has a perfect pairing on cohomology.
 However, one can easily show that Poincar\'e duality on chain level implies Poincar\'e duality duality on cohomology (see~\cite[Lemma~2.4]{triplepaper}), so \emph{a dPD algebra is also an oriented PDGA}.

We will adhere to the following convention unless otherwise stated:

\begin{convention}
	We suppose that all pairings $\la-,-\ra$ are graded symmetric and have degree $n\in\N_0$.
	We use the correspondence~\eqref{Eq:pairing-orientation} of pairings and orientations without further noting it.
\end{convention}

We use the following notion of a Hodge decomposition; it was used in~\cite{Cieliebak2015} and \cite{Van2019} and under the name ``harmonious Hodge decomposition'' in \cite{LazarevHodge}.

\begin{definition}\label{Def:HodgeDecomposition}
Let $(V,\Dd,\la-,-\ra)$ be a cochain complex (or a DGA) with a pairing.
A \emph{Hodge decomposition} of $(V,\Dd,\la-,-\ra)$ is a direct sum decomposition $V = \HH\oplus\im\Dd\oplus C$, where~$\HH$ is a complement of $\im\Dd$ in $\ker\Dd$, $C$ a complement of $\ker\Dd$ in~$V$, and $C \perp C\oplus\HH$.
We call $\HH$ the \emph{harmonic subspace} and $C$ the \emph{coexact part} of the Hodge decomposition.
We say that $(V,\Dd,\la-,-\ra)$ is \emph{of Hodge type} if it admits a Hodge decomposition.
\end{definition}

\begin{example}\label{Ex:DeRham}
	The term ``Hodge decomposition'' usually refers to the $L_2$-orthogonal decomposition $\Om = \HH\oplus \im\Dd\oplus \im \Dd^*$ of the graded vector space of smooth \emph{de Rham forms} $\Om\coloneqq\Om(X)$ on an oriented closed $n$-manifold~$X$ equipped with a \emph{Riemannian metric}~$g$.
	The \emph{$L_2$-inner product} is for all $\omega_1, \omega_2\in\Om$ defined by $(\omega_1,\omega_2)=\int_X \omega_1\wedge\mathop{\star}\omega_2$, where $\star=\star_g$ is the \emph{Hodge star} and $\wedge$ the \emph{wedge product}.
	Further, $\Dd\colon\Om \to \Om$ is the \emph{de Rham differential}, $\Dd^*\colon \Om\to \Om$ the \emph{codifferential}, and $\HH=\ker \Dd \cap \ker \Dd^*$ is the space of \emph{harmonic forms}, which is finite-dimensional and isomorphic to the \emph{de Rham cohomology} $\HDR(X)\coloneqq \Hh(\Om(X),\Dd)$ (see \cite[\S 6.8]{Warner1983}).
	Recall that $\mathop{\star}\mathop{\star}=\pm\Id$ and $\Dd^*=\pm\mathop{\star}\Dd\mathop{\star}$, where the sign depends on the degree.

	The \emph{de Rham complex} $(\Om,\Dd)$ equipped with the wedge product $\wedge$ is a unital CDGA over~$\R$ called the \emph{de Rham algebra}.
	\emph{Stokes' theorem} (see \cite[\S 4.9]{Warner1983}) implies that the \emph{integration} $\Or\coloneqq\Or_X\colon\Om\to \R$, $\omega\mapsto \int_X \omega$ defines an orientation on $(\Om,\Dd)$.
	We refer to the oriented CDGA $(\Om,\Dd,\Or)$ as \emph{oriented de Rham algebra} and to the corresponding pairing $\la\omega_1,\omega_2\ra =\int_X \omega_1\wedge \omega_2$ as \emph{intersection pairing}.  
	One can easily check that $\Om=\HH\oplus \im\Dd\oplus \im \Dd^*$ is a Hodge decomposition of $(\Om,\Dd,\la-,-\ra)$ in the sense of Definition~\ref{Def:HodgeDecomposition}, so $(\Om,\Dd,\Or)$ is of Hodge type.
	Moreover, $\la-,-\ra\colon \Om\times\Om\to\R$ is nondegenerate, and, using the Hodge decomposition together with $\HH\simeq \HDR\coloneqq \HDR(X)$ and $\dim\HH<\infty$, we deduce that $\la -,-\ra_*\colon \HDR\times \HDR\to\R$ is perfect.
	To sum up, the oriented de Rham algebra $(\Om(X),\Dd,\Or_X)$ of an oriented closed $n$-manifold~$X$ is an \emph{oriented PDGA of Hodge type of degree $n$ with a nondegenerate pairing}.
\end{example}

A Hodge decomposition can also be constructed using \emph{Poincar\'e duality}: 

\begin{lemma}[{\cite[Lemma~11.1]{Cieliebak2015},\cite{triplepaper}}]\label{Lem:CyclicIsHodge}
	If $\Char(\K)\neq 2$, then any cyclic cochain complex (or a DGA) $(V,\Dd,\la-,-\ra)$ is of Hodge type.\hfill\qed
\end{lemma}

\begin{remark}
	The proof of \cite[Lemma~11.1]{Cieliebak2015} starts by constructing a Hodge star operator $\star$ on $V$ with respect to an auxiliary inner product $(-,-)$ and continues by standard arguments for the Hodge decomposition from Example~\ref{Ex:DeRham}.
	The fact that $(-,-)$ is positive-definite is crucial, so the proof is applicable only when $\K$ is an \emph{ordered field}.
	In~\cite{triplepaper} we give an alternative proof based on the construction of a special propagator $\Pp$ (see Definition~\ref{Def:Propagator}).
	This proof only requires that $\Char(\K)\neq 2$ (due to the symmetrization $\Pp$).
	In Example~\ref{Ex:NoHodge} we give a counterexample to Lemma~\ref{Lem:CyclicIsHodge} for $\Char(\K)=2$ and even $n\in \N$.
\end{remark}
\begin{remark}\label{rem:existencehodge}
	Suppose that $(V,\Dd,\la-,-\ra)$ is a cochain complex (or a DGA) with a pairing of Hodge type such that the induced pairing $\la-,-\ra_*\colon \Hh(V)\times \Hh(V)\to \K$ is perfect (e.g., $V$ is a PDGA of Hodge type).
	Suppose that $\Char(\K)\neq 2$.
	Then \cite[Remark~2.6]{Van2019} shows that for any complement $\HH$ of $\im\Dd$ in $\ker\Dd$ there is a Hodge decomposition $V=\HH\oplus\im\Dd\oplus C$. 
\end{remark}

We differ from \cite{Lambrechts2007} by using the following notion of morphism of PDGAs:

\begin{definition}
	Given PDGAs $(V,\Dd,\Or_*)$ and $(V^\prime,\Dd^\prime,\Or^\prime_*)$ of the same degree $n\in\N_0$, a \emph{PDGA morphism} $f\colon V\to V^\prime$ is a DGA morphism such that the induced map on cohomology $f_*\colon \Hh(V)\to\Hh(V^\prime)$ satisfies $\Or_*^\prime\circ f_* = \Or_*$.
\end{definition}

\begin{remark}\label{Rem:Degree}
Suppose that $(V,\Dd,\Or_*)$ and $(V^\prime,\Dd^\prime,\Or_*^\prime)$ are PDGAs of degree $n\in\N_0$ with connected cohomology.
Poincar\'e duality then gives unique generators $h\in \Hh^n(V)$ and $h^\prime \in \Hh^n(V^\prime)$ such that $\Or_*(h) =\Or_{*}^\prime(h^\prime) = 1$.
If $f\colon (V,\Dd)\to (V^\prime,\Dd^\prime)$ is any DGA morphism, then we have $f_*(h) = \lambda h^\prime$ for some $\lambda=\lambda(f)\in\K$, which we refer to as \emph{mapping degree} of $f$.
We see that $f$ is a PDGA morphism if and only if $\lambda=1$.
Note that we have $\la f_*(h_1),f_*(h_2)\ra_*^\prime = \lambda \la h_1, h_2\ra_*$ for all $h_1,h_2\in \Hh(V)$.

In geometric applications, if $f\colon X \to X^\prime$ is a smooth homotopy equivalence of connected oriented closed $n$-manifolds, then the pullback $f^*\colon \Om(X^\prime)\to \Om(X)$ satisfies $\lambda(f^*)\in \{\pm 1\}$ (proof using the cohomology with $\Z$-coefficients).
\end{remark}

A cohomology orientation $\Or_*\colon \Hh(V)\to\K$ of degree $n$ can always be extended to a chain-level orientation $\Or\colon V\to \K$ of degree $n$ such that $\Or_*([v])=\Or(v)$ for all $[v]\in\Hh(V)$.
Such extension is unique if and only if $V^{n}\subset\ker\Dd$.
This occurs, for instance, when $V$ is a dPD algebra of degree $n$.
We then have the following:

\begin{lemma}\label{Lem:dPDMorphism}
	Let $(V,\Dd,\Or)$ and $(V^\prime,\Dd^\prime,\Or^\prime)$ be oriented DGAs of degree $n\in\N_0$ such that the pairing $\la-,-\ra\colon V\times V\to \K$ is nondegenerate. 
	Suppose that $f\colon (V,\Dd)\to (V^\prime,\Dd^\prime)$ is a DGA morphism that preserves cohomology orientation.
	Then $f$ preserves chain-level orientation and is injective.
\end{lemma}

\begin{proof}
Because $V$ is non-negatively graded and $\la-,-\ra\colon V^k \times V^{n-k}\to \K$ is nondegenerate, we have $V^k = 0$ whenever $k\not\in[0,n]$.
Given $v\in V^n$, we have $\Dd v \in V^{n+1}= 0$, hence $\Dd f(v) = f(\Dd v) = 0$, so
\begin{equation}\label{eq:OrientationInTopDegree}
	\Or^\prime(f(v)) = \Or^\prime_*([f(v)]) = \Or^\prime_*(f_*([v])) = \Or_*([v]) = \Or(v).
\end{equation}
Therefore, $f$ preserves chain-level orientation.
If $v_1\in V^k$ is such that $f(v_1)=0$, then for all $v_2\in V^{n-k}$ we have 
\[
\la v_1,v_2\ra = \Or(v_1 v_2) = \Or^\prime(f(v_1 v_2)) = \Or^\prime(f(v_1)f(v_2)) = 0,
\]
so the nondegeneracy of $\la-,-\ra$ implies $v_1=0$.
Therefore, $f$ is injective.
\end{proof}

The following notions from \cite{Van2019} are crucial for our construction:

\begin{definition}\label{Def:NonDegQuotient}
Let $(V,\Dd,\la-,-\ra)$ be a cochain complex (or a DGA) with a pairing.
We define the \emph{degenerate subspace} $V_\rmdeg$, resp., the \emph{nondegenerate quotient} $\QQ(V)$ by
\begin{equation*}
	V_\rmdeg \coloneqq \{v\in V \mid v \perp V\},\quad\text{resp.,}\quad\QQ(V) \coloneqq V / V_\rmdeg.
\end{equation*}
We denote by $\pi_\QQ\colon V\onto \QQ\coloneqq\QQ(V)$ the \emph{quotient map} and by $\la-,-\ra_{\QQ}\colon \QQ\times\QQ\to\K$ the unique pairing such that $\la\pi_\QQ(v_1),\pi_\QQ(v_2)\ra_\QQ=\la v_1,v_2\ra$ for all $v_1, v_2\in V$.
\end{definition}

Because $V_\rmdeg\subset V$ is a subcomplex (or a dg-ideal), there is a unique structure of a cochain complex (or a DGA) on $\QQ(V)$ such that $\pi_\QQ\colon V\onto\QQ(V)$ is a chain map (or a DGA morphism).
Statement~(a) of the following Lemma corresponds to \cite[Lemma~2.8]{Van2019}.
Statement (b) can be seen as a consequence of~\cite[Theorem~2.7]{LazarevHodge}, which asserts that~$V$ is of Hodge type if and only if $V_\rmdeg$ admits a finite-dimensional dg-complement.

\begin{lemma}\label{Lem:HodgeType}
 	Let $(V,\Dd,\la-,-\ra)$ be a cochain complex (or a DGA) with a pairing.
	Then the following holds:
	\begin{enumerate}[label=(\alph*)]
		\item If $(V,\Dd,\la-,-\ra)$ is of Hodge type and $\la -,-\ra_*\colon\Hh(V)\times\Hh(V)\to\K$ is nondegenerate, then $\pi_\QQ\colon V\onto\QQ(V)$ is a quasi-isomorphism.
		\item Suppose that $\Char(\K)\neq 2$.
			If $\pi_\QQ\colon V\onto\QQ(V)$ is a quasi-isomorphism and $\QQ(V)$ is of finite type, then $(V,\Dd,\la-,-\ra)$ is of Hodge type.
	\end{enumerate}
\end{lemma}

\begin{proof}

As for (a), consider a Hodge decomposition $V=\HH\oplus\im\Dd\oplus C$.
Let $v\in V_\rmdeg$ be such that $\Dd v=0$.
The relation $v\perp V$ implies $[v]\perp \Hh(V)$, so the nondegeneracy of $\la-,-\ra_*$ yields $[v]=0$ in $\Hh(V)$.
Therefore, there exists $c\in C$ such that $v=\Dd c$.
For any $h\in \HH$ and $c^\prime, c^\dprime\in C$ we use $C\perp C\oplus \HH$ to compute
\begin{equation*}
\la h + \Dd c^\prime + c^\dprime, c \ra = \la \Dd c^\prime, c \ra = \pm \la c^\prime, \Dd c\ra = \pm \la c^\prime, v \ra = 0.
\end{equation*}
This implies that $c\in V_\rmdeg$, so $[v]=0$ in $\Hh(V_\rmdeg)$.
We conclude that $\Hh(V_\rmdeg) = 0$.
However, $\pi_\QQ\colon V\onto\QQ(V)$ is surjective and $V_\rmdeg=\ker \pi_\QQ$, so $\Hh(V_\rmdeg)=0$ if and only if $\pi_\QQ$ is a quasi-isomorphism.

As for (b), we assume that $\Hh(V_\rmdeg)=0$, so there is a subcomplex $W\subset V$ such that~$V$ can be written as direct sum of cochain complexes $V=V_\rmdeg\oplus W$ (we say that $W$ is a dg-complement of $V_\rmdeg$).%
\footnote{If we work over a field $\K$, then every acyclic subcomplex admits a dg-complement (see \cite[Lemma~6.1.12]{MyPhD} for an elementary proof).} 
The restriction of $\pi_\QQ\colon V\to \QQ(V)$ to~$W$ is an isomorphism of cochain complexes with pairings.
Because $\la-,-\ra_\QQ$ is nondegenerate and $\QQ(V)$ is of finite type, $\la-,-\ra_\QQ$ is perfect, so $\la-,-\ra|_{W\times W}$ is perfect.
Hence, $(W,\Dd|_W,\la-,-\ra|_{W\times W})$ is a cyclic cochain complex, and Lemma~\ref{Lem:CyclicIsHodge} gives a Hodge decomposition $W=\HH\oplus\im\Dd^\prime \oplus C^\prime$.
Let $C^{\prime\prime}\subset V_\rmdeg$ be any complement of $\Dd V_\rmdeg\subset V_\rmdeg$.
Setting $C\coloneqq C^\prime \oplus C^{\prime\prime}$, we have $C\perp C\oplus \HH$, so $V=\HH\oplus\im\Dd\oplus C$ is a Hodge decomposition.
\end{proof}

\begin{remark}
In \cite[Proposition 6.1.17]{MyPhD} we show that a Hodge decomposition $V=\HH\oplus\im\Dd\oplus C$ induces a Hodge decomposition $\QQ(V)=\pi_\QQ(\HH)\oplus\im\Dd_\QQ \oplus \pi_\QQ(C)$.
\end{remark}

\begin{definition}\label{Def:Propagator}
	Let $(V,\Dd,\la-,-\ra)$ be a cochain complex (or a DGA) with a pairing.

	A \emph{special propagator} is a linear map $\Pp\colon V\to V$ with $\deg\Pp=-1$ that satisfies the relations
	\begin{equation}\label{Eq:SpecialPropagator-rel}
		\Dd\circ\Pp\circ\Dd = -\Dd,\quad \Pp\circ\Dd\circ\Pp = -\Pp,\quad\text{and}\quad\Pp\circ\Pp = 0,
	\end{equation}
	together with the following \emph{symmetry condition} for all $v_1, v_2\in V$:
	\begin{equation}\label{Eq:self-adjoint}
		\la \Pp v_1, v_2 \ra = (-1)^{\deg v_1} \la v_1, \Pp v_2 \ra.
	\end{equation}
	Without the symmetry condition, $\Pp$ would be called a \emph{special homotopy operator}.

	A \emph{harmonic projection} is a \emph{linear projection} $\pi\colon V\to V$---that is, we have $\deg\pi=0$ and $\pi\circ\pi=\pi$---that satisfies the relations 
	\[
		\ker \Dd = \im \pi \oplus \im \Dd\quad\text{and}\quad \im\Dd \subset \ker\pi,
	\]
	together with the following \emph{symmetry condition} for all $v_1, v_2\in V$:
	\begin{equation*}
		\la \pi(v_1),v_2\ra = \la v_1,\pi(v_2)\ra.
	\end{equation*}
\end{definition}

Given a special propagator $\Pp\colon V\to V$, it is straightforward to check that the linear map $\pi_\Pp\colon V\to V$ defined by
\begin{equation}\label{Eq:AssociatedProjection}
	\pi_\Pp\coloneqq \Id + \Dd \circ \Pp + \Pp \circ \Dd
\end{equation}
is a harmonic projection and that the following relations hold:
\[
	\Pp\circ\pi_\Pp = \pi_\Pp\circ\Pp = 0.
\]
The proof of the following lemma is also straightforward:

\begin{lemma}\label{Lem:HodgeDecompositionsAndPropagators}
Let $(V,\Dd,\la-,-\ra)$ be a cochain complex (or a DGA) with a pairing.
To each Hodge decomposition $V=\HH\oplus\im\Dd\oplus C$ we can associate a special propagator $\Pp\colon V\to V$ defined for all $v\in V$ as follows:
\begin{equation}\label{Eq:SpecialPropagator}
	\Pp(v) \coloneqq \begin{cases} -c, & \text{if }v=\Dd c\text{ for some }c\in C,\\
	0, & \text{for }v\in \HH\oplus C.
	\end{cases}
\end{equation}
Conversely, to each special propagator $\Pp\colon V\to V$ we can associate a Hodge decomposition $V = \HH_\Pp \oplus \im\Dd\oplus C_\Pp$ given by
\begin{equation}\label{Eq:Correspondence}
	\HH_\Pp \coloneqq \im\pi_\Pp\quad\text{and}\quad C_\Pp \coloneqq \im\Pp.
\end{equation}
This gives a one-to-one correspondence of Hodge decompositions of $(V,\Dd,\la-,-\ra)$ and special propagators in $(V,\Dd,\la-,-\ra)$.
\qed
\end{lemma}

The next concept is of prime importance in~\cite{Van2019}:

\begin{definition}\label{Def:SmallSubalgebra}
Let $(V,\Dd,\la-,-\ra)$ be a DGA with a pairing of Hodge type.
Let $\Pp\colon V\to V$ be a special propagator, and let $\HH_\Pp = \im\pi_\Pp$ be the harmonic subspace of the associated Hodge decomposition.
The \emph{small subalgebra} $\SS_\Pp\subset V$ is defined as the least dg-subalgebra $\SS \subset V$ such that
\begin{equation}\label{Eq:SmallSubalgebra}
	\HH_\Pp \subset \SS\quad\text{and}\quad \Pp(\SS)\subset \SS.
\end{equation}
\end{definition}

\begin{remark}
Clearly, such $\SS_\Pp$ exists and is uniquely given as the intersection of all dg-subalgebras $\SS\subset V$ that satisfy \eqref{Eq:SmallSubalgebra}.
\end{remark}

 Claims (a) and (b) of the following lemma essentially correspond to \cite[Proposition~3.3]{Van2019}, where a degreewise description of $\SS_\Pp$ for a 1-connected PDGA $V$ is provided.
 We re-establish these statements by identifying~$\SS_\Pp$ with the colimit of a filtration induced by decorated binary trees known from the \emph{A$_\infty$-homotopy transfer}. 

\begin{lemma}\label{Lem:SmallSubalgebra}
Let $(V,\Dd,\la-,-\ra)$ be a DGA with a pairing of Hodge type, and let $\Pp\colon V\to V$ be a special propagator.
Then the following holds:
\begin{enumerate}[label=(\alph*)]
\item As a graded vector space, $\SS_\Pp$ is generated by iterated applications of $\Pp$ and the product~$\cdot$~to $k$-tuples of elements $h_1,\dotsc,h_k\in \HH_\Pp$ for $k\in \N$.
\item Suppose that $\Hh(V)$ is 1-connected.
      Then~$\SS_\Pp$ is also 1-connected; if $\Hh(V)$ is additionally of finite type, then so is~$\SS_\Pp$.
\item The restriction $\Pp\vert_{\SS_\Pp}\colon \SS_\Pp\to\SS_\Pp$ is a special propagator in the DGA with a paring $(\SS_\Pp, \Dd\vert_{\SS_\Pp},\la-,-\ra\vert_{\SS_\Pp})$, and the associated Hodge decomposition reads $\SS_\Pp= \HH_\Pp \oplus \Dd \SS_\Pp \oplus \Pp(\SS_\Pp)$.
In particular, the inclusion $\iota\colon \SS_\Pp \into V$ is a pairing preserving DGA quasi-isomorphism.
\end{enumerate}
\end{lemma}

\begin{proof}
We denote by $\Trees_k$ the set of isotopy classes of embeddings into~$\R^2$ of rooted binary trees with $k\in \N$ leaves together with a labeling of each edge (both interior and leaf edges) by either~$\Pp$ or the identity~$\Id$ and a labeling of each interior vertex by the product~$\cdot$~.
We orient each edge of $T\in\Trees_k$ toward the root and order its leaves in the counterclockwise direction in $\R^2$, so that we can interpret $T$ as a composition rule for the operations $\Id$, $\Pp$, and~$\cdot$~.
We can thus assign to $T$ a linear map $\Ev_{T}\colon V^{\otimes k}\to V$ in a natural way.
We define a filtration $(\Filtr_k)_{k\in\N}$ on $\SS_\Pp$ by
\begin{equation*}
	\Filtr_k \coloneqq \sum_{i=1}^k \sum_{T\in \Trees_i} \Ev_T(\HH^{\otimes i}_\Pp).
\end{equation*}
We see that $\Filtr_k\cdot \Filtr_{k^\prime}\subset \Filtr_{k+k^\prime}$ for all $k, k^\prime \in \N$ by grafting two trees to a common root.
We also see that $\Dd(\Filtr_k) \subset \Filtr_k$ by propagating $\Dd$ from the root to the leaves using the relation $\Dd \circ \cdot = \cdot \circ (\Dd \otimes\Id)+\cdot\circ(\Id\otimes\Dd)$ at an interior vertex and the relation $\Dd\circ \Pp = \pi_\Pp - \Pp\circ\Dd - \Id$ at an edge labeled by~$\Pp$.
Finally, we see that $\Pp(\Filtr_k) \subset \Filtr_k$ because the root edge is labeled either by~$\Pp$, so that another application of~$\Pp$ gives~$0$ since $\Pp \circ \Pp = 0$, or by~$\Id$, so that another application of~$\Pp$ leads to its relabeling by~$\Pp$. 
Therefore, the colimit $\bigcup_{k=1}^\infty \Filtr_k$ is a dg-subalgebra of $(V,\Dd)$ that contains $\HH_\Pp$ and is closed under~$\Pp$.
The minimality of~$\SS_\Pp$ implies that $\SS_\Pp = \bigcup_{k=1}^\infty \Filtr_k$, i.e., that $(\Filtr_k)_{k\in\N}$ is exhaustive.
Claim~(a) follows because~$\Filtr_k$ is clearly generated by iterated applications of $\Pp$ and~$\cdot$~to tuples of elements of~$\HH_\Pp$ of length~$\le k$.

For $k\ge 2$, consider an element of the quotient $\Filtr_k/\Filtr_{k-1}$ that is represented by $\Ev_T(h_1,\dotsc,h_k)\in \Filtr_k$ for some $T\in\Trees_k$ and $h_1,\dotsc,h_k\in \HH_\Pp$.
We can suppose that the leaf edges of $T$ are labeled by~$\Id$ since otherwise we would have $\Ev_T(\HH^{\otimes k}_\Pp)=0$ using $\Pp\circ\pi_\Pp=0$.
Suppose that the leaves~$i$ and $i+1$ are connected to the same interior vertex. 
We can then write 
\[
	\Ev_T(h_1,\dotsc,h_k) = \Ev_{T^\prime}(h_1, \dotsc, h_{i-1}, h_i \cdot h_{i+1},h_{i+2}, \dotsc, h_k)
\]
for the tree $T^\prime\in \Trees_{k-1}$ obtained by combining the two leaves into one.
If $h_i\cdot h_{i+1}\in\HH_\Pp$, then $\Ev_T(h_1,\dotsc,h_k)\in \Filtr_{k-1}$, so the corresponding class in $\Filtr_k/\Filtr_{k-1}$ is trivial.
This is clearly the case if $h_i$ or~$h_{i+1}$ is a multiple of $1\in V$. 
Together with the assumption that $\HH_\Pp\simeq \Hh(V)$ is 1-connected this means that we can assume that $\deg h_i \ge 2$ for all $i\in \{1,\dotsc,k\}$.
Further, any rooted binary tree with~$k$ leaves has~$2k-1$ edges out of which $k$ are leaf edges.
Because $\deg \Pp = -1$, we obtain the following degree estimate:
\begin{equation}\label{Eq:Estimate}
		\deg\bigl(\Ev_{T}(h_1,\dotsc,h_k)\bigr) \ge \sum_{i=1}^k \deg h_i - (k-1)\ge  k + 1.
\end{equation}
This together with $\Filtr_1 = \HH_\Pp$ implies that $\SS_\Pp$ is 1-connected. 
If $\Hh(V)\simeq \HH_\Pp$ is of finite type, then~$\Filtr_k$ is of finite type for each $k\in\N$ because~$\Trees_k$ is finite.
Now, \eqref{Eq:Estimate} implies that the filtration $(\Filtr_k^i)_{k\in\N}$ stabilizes in each degree $i\in \N_0$, so $\SS_\Pp = \bigcup_{k=1}^\infty \Filtr_k$ is of finite type, too.
Claim (b) is thereby proven.

As for (c), the restriction $\Pp|_{\SS_\Pp}$ is again a special propagator since $\SS_\Pp\subset V$ is an oriented dg-subalgebra.
The corresponding Hodge decomposition is given by~\eqref{Eq:Correspondence}.
\end{proof}

\section{Extension of Hodge type}\label{Sec:4}
This section contains the main construction of this paper, Lemma~\ref{Lem:Extension}, and its consequences: Proposition~\ref{Prop:ExtensionSimplyConnected} in the 1-connected case and Proposition~\ref{Prop:ExtensionConnected} in the connected case of finite type.
We also introduce the notions of pre-Hodge decomposition, Hodge twist, and PDGA retraction to streamline the presentation.

\begin{definition}\label{Def:PreHodge} 
A \emph{pre-Hodge decomposition} of a cochain complex (or a DGA) with a pairing $(V,\Dd,\la-,-\ra)$ is a direct sum decomposition $V=\HH\oplus \im\Dd \oplus C$, where~$\HH$ is a complement of $\im\Dd$ in $\ker\Dd$ and $C$ is a complement of $\ker\Dd$ in $V$ such that $\HH\perp C$.
We call $\HH$ the \emph{harmonic subspace} and $C$ the \emph{coexact part} of the pre-Hodge decomposition.
\end{definition}

\begin{remark}\label{Rem:PreHodge}
	One can easily see that the data of a pre-Hodge decomposition $V=\HH\oplus\im\Dd\oplus C$ is equivalent to the data of a harmonic projection $\pi\colon V \onto\HH\subset V$ together with a choice of a complement $C$ of $\im \Dd$ in $\ker\pi=\HH^\perp$.
	In analogy with Lemma~\ref{Lem:HodgeDecompositionsAndPropagators}, this is equivalent to the data of a special homotopy operator $\Pp\colon V\to V$ such that the map $\pi_\Pp\colon V\to V$ associated to $\Pp$ in~\eqref{Eq:AssociatedProjection} is a harmonic projection. 
\end{remark}

\begin{lemma}\label{Lem:PreHodge}
	Let $(V,\Dd,\la-,-\ra)$ be a cochain complex (or a DGA) with a pairing such that the induced pairing on cohomology $\la-,- \ra_*\colon \Hh(V)\times\Hh(V)\to \K$ is perfect.
	Then each complement $\HH$ of $\im \Dd$ in $\ker\Dd$ admits a unique harmonic projection $\pi_\HH\colon V\onto\HH\subset V$.
	In particular, there exists a pre-Hodge decomposition $V=\HH\oplus\im\Dd\oplus C$.
\end{lemma}
\begin{proof}
	One can easily see that $\pi_\HH$ must be the orthogonal projection onto $\HH$, which exists and is unique because the restriction $\la-,-\ra|_{\HH\times \HH}$ is perfect.
\end{proof}

In order to motivate our method, let us think how to modify (or ``twist'') the coexact part $C$ of a pre-Hodge decomposition $V=\HH\oplus\im\Dd\oplus C$ in order to obtain a Hodge decomposition $V=\HH\oplus\im\Dd\oplus C^\prime$. 
Inspired by \cite[Remark~2.6]{Van2019}, we look for~$C^\prime$ in the form 
\begin{equation}\label{Eq:graph}
C^\prime=\mathrm{graph}(\mu)\coloneqq\{ c + \mu(c) \mid c\in C\}
\end{equation}
for a linear map $\mu\colon C\to\im\Dd$ of degree~$0$.
Using $C\perp\HH$ and $\im\mu\subset\im\Dd\perp\HH$ we get $C^\prime\perp \HH$.
Moreover, using $\im\Dd\perp\im\Dd$ we get that $C^\prime\perp C^\prime$ holds if and only if $\mu$ satisfies the following equation for all $c_1, c_2\in C$:
\begin{equation}\label{Eq:HodgeTwistGeneral}
	\la \mu(c_1),c_2 \ra + \la c_1,\mu(c_2)\ra + \la c_1, c_2\ra = 0.
\end{equation}
In fact, if $\la-,-\ra_*\colon \Hh(V)\times\Hh(V)\to \K$ is perfect, then $\im\Dd\oplus C = \HH^\perp = \im\Dd\oplus C^\prime$, so~\eqref{Eq:graph} gives a one-to-one correspondence of Hodge decompositions $V=\HH\oplus\im\Dd\oplus C^\prime$ and solutions $\mu\colon C\to \im\Dd$ of \eqref{Eq:HodgeTwistGeneral}.

Clearly, a pre-Hodge decomposition $V=\HH\oplus\im\Dd\oplus C$ is a Hodge decomposition if and only if the following nondegenerate quotient vanishes:
\[
	\QQ(C) \coloneqq C/C_\rmdeg,\quad\text{where}\quad C_\rmdeg \coloneqq \{c\in C \mid c\perp C\}.
\]
The pairing $\la-,-\ra_\QQ\colon \QQ(C)\times \QQ(C)\to \K$ is nondegenerate, so we have $\QQ(C) = \QQ^0(C)\oplus \dotsb\oplus \QQ^n(C)$.
If in addition $\dim \QQ(C)<\infty$, then \emph{Poincar\'e duality} holds: 
\[ 
	\QQ^i(C)\simeq \QQ^{n-i}(C).
\]
This suggests to restrict to solutions $\mu\colon C\to \im\Dd$ of~\eqref{Eq:HodgeTwistGeneral} that vanish on~$C_\rmdeg$ and~$C^{[0,n/2)}$.
This leads to the following definition:

\begin{definition}\label{Def:HodgeTwist}
	Let $(V,\Dd,\la-,-\ra)$ be a cochain complex (or a DGA) with a pairing of degree $n\in\N_0$, and let $V=\HH\oplus\im\Dd\oplus C$ be a pre-Hodge decomposition.
	A~\emph{Hodge twist} of degree $k\in[n/2,n]$ is a linear map $\mu\colon\QQ(C)\to \im\Dd$ of degree~$0$ that satisfies the following conditions:
	\begin{enumerate}
		\item For all $i\in (n/2,k]$ and $c_1\in C^i$, $c_2\in C^{n-i}$ we have
		\begin{equation*}
			\la \mu(c_1), c_2 \ra + \la c_1, c_2\ra = 0.
		\end{equation*}
	\item If $n$ is even, then for all $c_1,c_2\in C^{n/2}$ we have
		\begin{equation*}	
			\la \mu(c_1),c_2\ra + \la c_1, \mu(c_2)\ra + \la c_1, c_2 \ra = 0.
		\end{equation*}
	\item For all $i\not\in[n/2,k]$ we have $\mu=0$ on $\QQ^{i}(C)$.
	\end{enumerate}
	Here we implicitly evaluate $\mu$ on $C$ via the quotient map $\pi_\QQ\colon V\to \QQ(V)$.
\end{definition}

A Hodge twist $\mu\colon \QQ(C)\to \im \Dd$ of degree $k$ can be viewed as a recipe to turn a pre-Hodge decomposition into a Hodge decomposition in degrees $i\in [n-k,k]$:

\begin{lemma}\label{Lem:HodgeTwist}
	Let $(V,\Dd,\la-,-\ra)$ be a cochain complex with a pairing of degree $n\in \N_0$.
	Suppose that $V=\HH\oplus\im\Dd\oplus C$ is a pre-Hodge decomposition and $\mu\colon \QQ(C)\to\im\Dd$ a Hodge twist of degree $k\in[n/2,n]$.
 	Then $V=\HH\oplus\im\Dd\oplus\mathrm{graph}(\mu)$ is a pre-Hodge decomposition which is Hodge in degrees $i\in [n-k,k]$.
\end{lemma}

\begin{proof}
	We saw that $V=\HH\oplus\im\Dd\oplus\mathrm{graph}(\mu)$ is a Hodge decomposition if and only if \eqref{Eq:HodgeTwistGeneral} holds for all $c_1\in C^i$, $c_2\in C^{n-i}$ and $i\in[0,n]$.
	Since $\mu$ vanishes on~$C^{[0,n/2)}$, equation~\eqref{Eq:HodgeTwistGeneral} reads $\la \mu(c_1),c_2\ra + \la c_1,c_2\ra = 0$ for $i<n/2$ and $\la c_1,\mu(c_2)\ra + \la c_1,c_2 \ra = 0$ for $i>n/2$.
	The first condition corresponds to Definition~\ref{Def:HodgeTwist} (1), and the second condition is due to the graded symmetry of $\la-,-\ra$ equivalent to the first.
\end{proof}

We will use the following CDGAs, which were used in~\cite{Lambrechts2007} as well:

\begin{definition}\label{Def:CharDef}
A \emph{free CDGA} on generators $x_\alpha$ with $\deg x_\alpha\ge 1$ ($\alpha\in I$) over a field $\K$ is the free graded commutative $\K$-algebra $\Lambda \coloneqq \Lambda_\K(x_\alpha,\Dd x_\alpha \mid \alpha\in I)$ together with the differential $\Dd\colon\Lambda\to\Lambda$ such that $\Dd(x_\alpha) = \Dd x_\alpha$ and $\Dd(\Dd x_\alpha)=0$ for all $\alpha\in I$.

A \emph{free acyclic CDGA} $\ALambda\coloneqq \ALambda_\K(x_\alpha\mid \alpha\in I)$ over $\K$ on the above generators $x_\alpha$ is defined as the free CDGA $\ALambda\coloneqq \Lambda$ if $\Char(\K) = 0$, and is defined as the CDGA
\[ 
	\ALambda\coloneqq\Lambda\otimes\Lambda_\K(x_{\alpha,i},y_{\alpha,i}\mid i\in \N, \alpha\in I, \deg x_{\alpha} \in 2\N)
\]
with $\deg x_{\alpha,i} = \deg y_{\alpha,i} - 1 = i p \deg x_{\alpha}$ and 
\[
	\Dd y_{\alpha,i} = x_{\alpha,i-1}^p, \quad \Dd x_{\alpha,i} = x_{\alpha,i-1}^{p-1}\Dd x_{\alpha,i-1}
\]
for all $i\in\N$ and $\alpha\in I$ if $\Char(\K)=p$.
Here we denote $x_{\alpha,0}\coloneqq x_{\alpha}$ and assume additionally that all elements of $\ALambda$ of odd degree square to $0$ if $p=2$.%
\footnote{This can be achieved by modding out the corresponding dg-ideal.
Note that if $x$ is a (non-exact) generator of odd degree, then $\Dd(x^2) = (\Dd x) x - x (\Dd x) = (\Dd x) x - (\Dd x) x = 0$, so $x^2$ represents a nontrivial cohomology class.}
\end{definition}

The following fact, which is used in~\cite{Lambrechts2007}, is proven here for the sake of completeness:

\begin{lemma}\label{Lem:FreeCDGA}
	Let $\K$ be a field, and let $x_\alpha$ with $\deg x_\alpha \ge 1$ ($\alpha\in I$) be a set of generators of a free CDGA $\Lambda\coloneqq \Lambda_\K(x_\alpha,\Dd x_\alpha\mid \alpha\in I)$.
	Let $\ALambda \coloneqq \ALambda_\K(x_\alpha\mid \alpha\in I)$ be the corresponding free acyclic CDGA.
	Then the following holds for the cohomology:
	\begin{enumerate}
		\item[(a)] $\Hh(\Lambda) = 0$ if $\Char(\K)=0$, and 
			\[ 
				H(\Lambda) = \Lambda_\K([x_\alpha^p],[x_\alpha^{p-1}\Dd x_\alpha]\mid \alpha\in I, \deg x_\alpha\in 2\N)
			\]
			if $\Char(\K)=p>2$. 
			This holds also for $p=2$ under the assumption that all elements of $\Lambda$ of odd degree square to $0$.
		\item[(b)] $\Hh(\ALambda) = \Span_\K\{[1]\}$.
\end{enumerate}
\end{lemma}

\begin{proof}
It suffices to prove the lemma for the case where $I=\{*\}$ is a singleton.
The general case can then be derived from K\"unneth's theorem. 
We denote $x\coloneqq x_{*}$ and $k \coloneqq \deg x + 1 \ge 2$ (we introduce $k$ to ensure compatibility with forthcoming notation).

Claim (a) is clear from the following computation:
\begin{align*}
	&k\text{ even:} &&\Lambda = \Span_\K\{x(\Dd x)^j, (\Dd x)^{j+1}\mid j\in\N_0\}, && \Dd(x (\Dd x)^{j})=(\Dd x)^{j+1},\\ 
	&k\text{ odd:} &&\Lambda = \Span_\K\{ x^j(\Dd x), x^{j+1}\mid j\in \N_0\}, &&\Dd(x^{j+1}) = (j+1) x^{j}\Dd x. 
\end{align*}

As for (b), we can assume that $k$ is odd, so we have $\ALambda = \Lambda_\K(x,\Dd x, x_i, y_i \mid i\in \N)$ with $\deg x_i = \deg y_i - 1 = ip(k-1)$ and $\Dd y_i = x_{i-1}^p$, $\Dd x_i = x_{i-1}^{p-1}\Dd x_{i-1}$ ($x_0\coloneqq x$).
For $i\in \N_0$ let $\ALambda_i\subset \ALambda$ be the dg-subalgebra generated by all $x_j, y_j$ with $j\le i$.
We set $y_0\coloneqq 0$, so $\ALambda_0 = \Lambda$.
We then have $\Hh(\ALambda_0)=\Lambda_\K([x^p],[x^{p-1}\Dd x])$ by (a).
We will show by induction that $\Hh(\ALambda_i)=\Lambda_\K([x_i^p],[x^{p-1}_i\Dd x_i])$ holds for all $i\in \N$. 
Given $i\in \N$ such that $\Hh(\ALambda_{i-1})=\Lambda_\K([x_{i-1}^p],[x^{p-1}_{i-1}\Dd x_{i-1}])$, consider for each $k\in\N_0$ the subcomplex $\ALambda_{i,k}\subset\ALambda_i$ spanned by words in generators $x_j$, $y_j$ ($j\le i$) such that the sum of powers of~$x_i$ and~$y_i$ does not exceed~$k$.
We obtain a filtration $\ALambda_{i,0}\subset \ALambda_{i,k}\subset\ALambda_{i,k+1}\subset\dotsb\subset \ALambda_i$.
Following \cite[Section~5]{Weibel}, the filtration induces a homological spectral sequence $(E_r,\partial_r)$ ($r\in\N_0$) such that~$E_0$ consists of elements of the graded module $\mathrm{gr}(\ALambda_i) \coloneqq \oplus_{k=1}^\infty (\ALambda_{i,k}/\ALambda_{i,k-1})$, and we have $E_{r+1}=\Hh(E_r,\partial_r)$, where $\partial_r$ is induced from $\Dd$ for each $r\in\N_0$ (because we are in the homological setting, the grading on~$\ALambda_i$ has to be reversed, so each~$E_r$ is supported in the fourth quadrant under the negative diagonal).
So, the first page $E_1$ consists of elements of
\[
	\Hh(\mathrm{gr}(\ALambda_i))\simeq \mathrm{gr}(\Lambda_\K(x_i,y_i))\otimes \Hh(\ALambda_{i-1}) \simeq\Lambda_\K(x_i,y_i)\otimes\Lambda_\K([x_{i-1}^p],[x_{i-1}^{p-1}\Dd x_{i-1}]),
\]
where $\Lambda_\K(x_i,y_i)$ is filtered by the total power of $x_i$ and $y_i$.
Together with the induced differential this is isomorphic to $\Lambda_\K(x_i,y_i,\Dd y_i, \Dd x_i)  \simeq \Lambda_\K(x_i,\Dd x_i)\otimes\Lambda_\K(y_i,\Dd y_i)$.
Therefore, the second page $E_2$ consists of elements of $\Hh(\Lambda_\K(x_i,y_i,\Dd x_i,\Dd y_i))=\Lambda_\K([x_i^{p}],[x_i^{p-1}\Dd x_i])$.
The induced differential clearly vanishes, so we have $\partial_2=0$, and the spectral sequence abuts at $E_2$.
Since the filtration is bounded (in each degree), \cite[Theorem~5.5.1]{Weibel} applies, and we conclude that $\Hh(\ALambda_i)=\Lambda_\K([x_i^{p}],[x_i^{p-1}\Dd x_i])$.
Because $\deg(x_i^p) = ip(k-1)\to \infty$ as $i\to \infty$, we have $\Hh^j(\ALambda) = \lim_{i\to\infty} \Hh^j(\ALambda_i) = 0$ for each $j\in \N$.
\end{proof}

Now comes our main construction, which is similar to the construction in~\cite{Lambrechts2007}. 

\begin{lemma}\label{Lem:Extension}
Let $(V,\Dd,\Or)$ be a 1-connected oriented PDGA of degree $n\ge 5$ over a field $\K$, let $V=\HH\oplus\im\Dd\oplus C$ be a pre-Hodge decomposition, and let $k\in [n/2,n]$.
For $k\ge n/2+1$ we additionally suppose that we are given a Hodge twist $\mu\colon\QQ(C)\to\im\Dd$ of degree $k-1$.
For \emph{even} $n$ we additionally assume the following conditions:
\begin{align}
	&\begin{gathered}[t]\text{If } k=n/2\text{ and }\Char(\K)=2,\text{ then } \dim \QQ^{n/2}(C)<\infty\\\text{and }\la q,q \ra_\QQ = 0\text{ for all } q\in \QQ^{n/2}(C).\label{Eq:SymplecticBasis}\tag{{$\ast$}}\end{gathered}\\
	&\text{If }k=n/2+1\text{, then }C^{n/2}\perp C^{n/2}. \label{Eq:AdditionalAssumption}\tag{{$\ast\ast$}}
\end{align}
Let $q_\alpha$ ($\alpha\in I$) be a basis of $\QQ^{k}(C)$, and let $\ALambda\coloneqq \ALambda_\K(x_\alpha\mid\alpha\in I)$ be a free acyclic CDGA with $\deg \Dd x_\alpha = \deg x_\alpha + 1 = k$. 
Consider the \emph{tensor product CDGA}%
\footnote{Our differential is the tensor product differential, whereas the differential in \cite{Lambrechts2007} is twisted.} 
\[
	(\wh V,\wh \Dd) \coloneqq (\ALambda,\Dd)\otimes (V,\Dd).
\]
We claim that $\wh V$ is 1-connected, that the canonical inclusion $\iota\colon V\into\wh V$ is a DGA quasi-isomorphism, and that there is an orientation $\wh \Or\colon \wh V \to \K$ of degree~$n$ such that $\wh \Or \circ \iota = \Or$.
We further claim that there is a pre-Hodge decomposition $\wh V = \HH\oplus\im\wh\Dd\oplus\wh C$ and a Hodge twist $\wh\mu\colon\QQ(\wh C)\to \im\wh\Dd$ of degree $k$ such that for $k\ge n/2+1$ the following \emph{extension property} holds for all $i\in[n/2,k-1]$:
\begin{equation}\label{eq:extension}
	\begin{gathered}
		\text{We have }C^i\subset \wh C^i\text{ and }C_\rmdeg^i \subset \wh C_\rmdeg^i\cap C^i,\\\text{the induced map }\QQ^{i}(C)\to \QQ^{i}(\wh C)
		\text{ is an isomorphism,}\\
	\text{ and we have }\wh\mu=\mu\text{ on }\QQ^{i}(\wh C)\simeq \QQ^{i}(C).	
	\end{gathered}
	\tag{$\dagger$}
\end{equation}
\end{lemma}

\begin{proof}
The fact that $\iota\colon V\into \wh V$ is a quasi-isomorphism follows from K\"unneth's theorem using that $\Hh(\ALambda)\simeq \K$ by Lemma~\ref{Lem:FreeCDGA}.
Using $n\ge 5$ we have 
\[
	\deg x_\alpha = k-1 \ge \lceil n/2\rceil -1 \ge 2,
\]
so~$\wh V$ is 1-connected.
Let us suppose that $\Char(\K)=0$ or $k = \deg x_\alpha + 1\in 2\N$, so $\ALambda = \Lambda\coloneqq \Lambda_\K(x_\alpha,\Dd x_\alpha)$.
We will address the case $\Lambda\subsetneq\ALambda$ at the end of the proof.

We start by defining the orientation $\wh \Or\colon \wh V\to \K$.
Having $\wh V = \Lambda \otimes V$, we set
\begin{equation*}
	\Lambda_{i,j}\coloneqq \Span_\K\bigl\{ x_{\alpha_1}\dotsm x_{\alpha_i} \Dd x_{\beta_1} \dotsm \Dd x_{\beta_j} \mid \alpha_1, \dotsc, \alpha_i, \beta_1, \dotsc, \beta_j\in I \bigr\}
\end{equation*}
for all $i,j \in \N_0$, and for all $i\in\N_0$ we define
\begin{equation*}
	\Lambda_i \coloneqq \bigoplus_{\substack{a+b=i}} \Lambda_{a,b}\quad\text{and}\quad\wh V_i \coloneqq \Lambda_i \otimes V.
\end{equation*}
Hence, we have $\wh V\simeq \oplus_{i\in\N_0}\wh V_i$ as cochain complexes.
We set $\wh \Or \coloneqq \Or$ on $\wh V_0 \simeq V$, $\wh \Or\coloneqq 0$ on $\wh V_i$ for each $i\ge 2$, and $\wh \Or \coloneqq 0$ on~$\wh V_1^j$ for each $j\neq n$.
In order to extend $\wh\Or$ to
\[
	\wh V_1^n = (\Lambda_{1,0}\otimes V^{n-(k-1)})\oplus(\Lambda_{0,1}\otimes V^{n-k}),
\]
we decompose each $v\in V$ as $v = h + \Dd c + c^\prime$ with unique $h\in \HH$ and $c, c^\prime \in C$ (recall that $V=\HH\oplus \im\Dd\oplus C$) and define for each $\alpha\in I$ the following:
\begin{equation}\label{Eq:DefOfOr}
	\begin{aligned}
		&\wh\Or(x_\alpha \cdot v) \coloneqq (-1)^k \Or(q_\alpha\cdot c)\quad\text{and}\\
		&\wh\Or(\Dd x_\alpha \cdot v ) \coloneqq \Or(q_\alpha \cdot c^\prime).
	\end{aligned}
\end{equation}
Here we identify each $q_\alpha\in \QQ^{k}(C) = C^k/C^k_\rmdeg$ with a representative $q_\alpha\in C^{k}$ for each $\alpha\in I$, and we will continue to do so further.
Using 
\[ \wh\Dd(\Dd x_\alpha \cdot v) = (-1)^k \Dd x_\alpha \cdot \Dd v = (-1)^k\wh\Dd(x_\alpha \cdot \Dd v)\]
one can see that
\begin{equation*}
	\wh \Dd \wh V_1 = \Span_\K \bigl\{ \wh\Dd(x_\alpha\cdot v) = \Dd x_\alpha\cdot v + (-1)^{k-1} x_\alpha\cdot \Dd v \mid \alpha\in I, v\in V \bigr\}.
\end{equation*}
Therefore, the equation $\wh\Or\circ\wh\Dd = 0$ holds on $\wh V_1$, and it evidently also holds on $\wh V_i$ for all $i\neq 1$.
We have $\wh\Or\circ\iota=\Or$ by construction, so $\wh\Or\colon\wh V\to\K$ is the required orientation.

We now construct the pre-Hodge decomposition $\wh V = \wh\HH\oplus\im\wh\Dd\oplus \wh C$.
We denote by $\la-,-\ra_{\wh \Or}\colon \wh V \times \wh V\to \K$ the pairing associated to $\wh \Or\colon \wh V\to \K$ via~\eqref{Eq:pairing-orientation} and by~$\perp_{\wh\Or}$ the corresponding perpendicularity relation.
Because the inclusion $\iota\colon V\into\wh V$ is a quasi-isomorphism, $\wh{\HH}\coloneqq\HH$ is a complement of $\im\wh \Dd$ in $\ker \wh \Dd$.
Lemma~\ref{Lem:PreHodge} applied to $(\wh V,\wh \Dd,\la-,-\ra_{\wh \Or})$ then guarantees the existence of a (unique) harmonic projection $\wh{\pi}_\HH\colon \wh V \onto \HH\subset \wh V$.
We denote 
\[
	\wh\pi_{\HH}^\perp\coloneqq \Id - \wh\pi_\HH\colon \wh V\to \wh V
\]
and set
\begin{equation}\label{eq:defofcoexact}
	\wh C_{\le 1} \coloneqq \wh\pi_{\HH}^\perp(C \oplus \Lambda_{1,0}\cdot V) = C\oplus \wh\pi_{\HH}^\perp(\Lambda_{1,0}\cdot V) \subset \HH \oplus C \oplus \Lambda_{1,0}\cdot V.
\end{equation}
Because $\ker \Dd \oplus C=V$, $\ker \wh \Dd_1 = \im \wh \Dd_1 = \wh\Dd(\Lambda_{1,0}\cdot V)$, and $\wh\Dd(\Lambda_{1,0}\cdot V)\oplus\Lambda_{1,0}\cdot V=\wh{V}_1$, we have that $C\oplus \Lambda_{1,0}\cdot V$ is a complement of $\ker\wh\Dd_{\le 1}=\ker \Dd \oplus \ker\wh\Dd_1$ in the subcomplex 
\[
	(\wh V_{\le 1} \coloneqq V\oplus \wh V_1,\wh\Dd_{\le 1}\coloneqq \Dd\oplus \wh \Dd_1)\subset (\wh V,\wh \Dd).
\]
Because $\im \wh\pi_\HH \subset \ker \wh \Dd_{\le 1}$, this holds also for $\wh{C}_{\le 1}$.
Further, for each $i\ge 2$ we pick an arbitrary complement~$\wh C_i$ of $\ker \wh \Dd_i$ in~$\wh V_i$, and we set
\begin{equation}\label{eq:defofcoexact2}
	\wh C \coloneqq \wh C_{\le 1} \oplus \bigoplus_{i=2}^\infty \wh C_i.
\end{equation}
We obtain a pre-Hodge decomposition $\wh V = \HH\oplus\im\wh\Dd\oplus \wh C$.%
\footnote{The application of $\wh\pi_{\HH}^\perp$ in the definition of $\wh C_{\le 1}$ seems to be necessary to ensure that $\wh C$ is perpendicular to $\HH$.
To illustrate this, let $v\in V$ and $h\in \HH$.
We have $\la x_\alpha \cdot v,h\ra_{\hat\Or} = \hat\Or(x_\alpha \cdot v\cdot h) = (-1)^{k} \Or(q_\alpha \cdot c)$, where we decomposed $v\cdot h = h^\prime + \Dd c + c^\prime$ for unique $h^\prime\in\HH$ and $c, c^\prime\in C$.
In order to get perpendicularity, we would need $c\in C_{\rmdeg}$, but we do not see a reason for that.}

We now construct the Hodge twist $\wh\mu\colon \QQ(\wh C)\to \im\wh\Dd$ and check~\eqref{eq:extension}.
We set $\wh \mu \coloneqq 0$ on $\QQ^{i}(\wh C)$ for all $i\not\in[n/2,k]$.
For $i\in[n/2,k]$, we have to consider the intermediate degrees $[n-k,n]$ of $\wh C$, which we spell out in Table~\ref{tab:coexact}.
\begin{table}
	\caption{The coexact part $\wh C$ in intermediate degrees $[n-k,k]$ provided that $V^0=\K$, $V^1=0$, and $n\ge 5$ (hence $k\ge 3$).}\label{tab:coexact}
	\begin{enumerate}
	\item[(a)] \underline{$n-(k-1)<k-1$}:
			\begin{align*}
				\wh C^i &= \begin{cases} 
		C^i\\ 
		C^i\\ 
		C^i\\ 
		C^{k-1} \oplus \Lambda_{1,0}\\ 
		C^k 
	\end{cases}&
					\wh C_\rmdeg^i &= \begin{cases}
		C^i_\rmdeg & i=n-k\\
		C^i_\rmdeg & i=n-(k-1)\\
		C^i_\rmdeg & n-(k-1)< i<k-1 \\ 
		C^{k-1}_\rmdeg \oplus \Lambda_{1,0} & i = k-1 \\
		C^k_{\rmdeg} & i=k
	\end{cases}
			\end{align*}
\item[(b)] \underline{$n=2(k-1)$}:
	\begin{align*}
	\wh C^{n-k} &= C^{n-k} &  \wh C^{n-k}_\rmdeg &= C_\rmdeg^{n-k} \\
	\wh C^{k-1} &=\wh C^{n-(k-1)}= C^{k-1}\oplus \Lambda_{1,0} &  \wh C^{k-1}_\rmdeg & = \wh C^{n-(k-1)}_\rmdeg = C_\rmdeg^{k-1}\oplus \Lambda_{1,0}\\
	\wh C^{k} &= C^{k} &  \wh C^k_\rmdeg &=C_\rmdeg^k
	\end{align*}
\item[(c)] \underline{$n=2k-1$}:
	\begin{align*}
		\wh C^{k-1} & = \wh C^{n-k} = C^{k-1} \oplus\Lambda_{1,0} & 
		\wh C^{k-1}_\rmdeg & = \wh C^{n-k}_\rmdeg = C_\rmdeg^{k-1} \oplus \Lambda_{1,0} \\
		\wh C^{k} & = \wh C^{n-(k-1)} = C^k & \wh C^{k}_\rmdeg & = \wh C^{n-(k-1)}_\rmdeg = C^k_\rmdeg 
	\end{align*}
\item[(d)] \underline{$n=2k$}:
	\begin{align*}
	  \wh C^{k}& = \wh C^{n-k} = C^{k} &  \wh C_\rmdeg^{k} &= \wh C^{n-k}_\rmdeg = C_\rmdeg^{k} 
	\end{align*}
\end{enumerate}
\end{table}
The table can be verified using~\eqref{eq:defofcoexact} and~\eqref{eq:defofcoexact2}.
Four cases (a)--(d) are distinguished.
We see immediately that in each (a)--(d) we have $\wh C^i_\rmdeg \cap C^i = C^i_\rmdeg$ and $\QQ^{i}(\wh C)\simeq \QQ^{i}(C)$ for all $i\in [n-k,k]$.
If $i\in [n/2,k)$, which occurs in (a) or (b), we consider the given Hodge twist $\mu\colon\QQ(C)\to \im\Dd$ and set $\wh\mu\coloneqq\mu$ on $\QQ^i(\wh C)$.
This definition assures that the extension property~\eqref{eq:extension} holds.
It remains to define~$\wh \mu$ on $\QQ^{k}(\wh C) \simeq \QQ^{k}(C)$ in each (a)--(d) and check that it is a Hodge twist of degree $k$.
In (a)--(c) we define for each $\alpha\in I$ the following:
\begin{equation}\label{Eq:DefOfMu}
	\wh \mu(q_\alpha) \coloneqq - \Dd x_\alpha.
\end{equation}
Suppose that we are in (d).
If $\Char(\K)\neq 2$, we define for each $\alpha\in I$ the following:
\begin{equation}\label{Eq:DefOfMuInTheMiddleDegree}
	\wh\mu(q_\alpha) \coloneqq -\frac{1}{2}\Dd x_\alpha.
\end{equation}
If $\Char(\K)=2$, we use~\eqref{Eq:SymplecticBasis} and assume that $q_\alpha$, $\alpha\in I = \{1,\dotsc,2l\}$ is a \emph{symplectic basis} of $\QQ^{n/2}(C)$ for some $l\in\N_0$---that is, $q_\alpha$ are perpendicular to each other except for $\la q_\alpha, q_{\alpha+l} \ra_\QQ = - \la q_{\alpha+l}, q_\alpha\ra_\QQ = 1$ for all $\alpha\in\{1,\dotsc,l\}$.%
\footnote{Recall that a \emph{symplectic vector space} $(\QQ,\la-,-\ra)$ is a vector space with a perfect antisymmetric pairing.
In order to construct a symplectic basis $q_1,\dotsc,q_{2l}$ of $(\QQ,\la-,-\ra)$, one picks some $q_1, q_{l+1}\in \QQ$ such that $\la q_1, q_{l+1} \ra = 1$ and considers the subspace $\QQ_1 \coloneqq \Span_\K\{q_1,q_{l+1}\}$.
The perfection of $\la-,-\ra$ implies $\dim \QQ_1 + \dim \QQ_1^\perp = \dim \QQ$, and the relations $\la q_1,q_{l+1} \ra= - \la q_{l+1},q_1\ra$ and $\la q_1,q_1\ra = \la q_{l+1},q_{l+1}\ra = 0$ imply $\QQ_1\cap \QQ_1^\perp = 0$.
If $\Char(\K)=2$, then $\la q,q\ra = 0$ does not automatically follow from the antisymmetry (which is the same as symmetry in this case), and has to be assumed independently.
Having $\QQ = \QQ_1\oplus \QQ_1^\perp$, the restriction of the pairing to $\QQ_1^\perp$ is perfect, and we can continue inductively with $\QQ_1^\perp$ instead of $\QQ$.}
We then replace~\eqref{Eq:DefOfMuInTheMiddleDegree} with the following definition for all $\alpha\in\{1,\dotsc,l\}$:
\begin{equation}\label{Eq:DefOfMuSympl}
	\wh\mu(q_\alpha) \coloneqq -\Dd x_{\alpha}\quad\text{and}\quad\wh\mu(q_{\alpha+l})\coloneqq 0.
\end{equation}
Having defined a linear map $\wh\mu\colon\QQ(\wh C)\to\im\wh\Dd$, we have to check that it is a Hodge twist of degree~$k$ according to Definition~\ref{Def:HodgeTwist}.

\begin{sublemma*}
For each $i\in [n/2,k]$ and all $c_1\in C^i$, $c_2\in C^{n-i}$ we have the following:
\begin{align}
	&\text{if }i>n/2: & 0=&\la\wh\mu(c_1),c_2\ra_{\wh\Or} + \la c_1,c_2\ra_{\wh\Or} \label{Eq:MidRel1},\\
	&\text{if }n\in 2\N\text{ and }i=n/2: & 0=&\la\wh\mu(c_1),c_2\ra_{\wh \Or} + \la c_1,\wh\mu(c_2)\ra_{\wh \Or} + \la c_1,c_2\ra_{\wh\Or}.\label{Eq:MidRel2}
\end{align}
\end{sublemma*}

\begin{proof}
First note that $\la c_1,c_2\ra_{\wh\Or}=\la c_1,c_2\ra$.
If $i<k$, then we have $\wh{\mu}=\mu$ on $\QQ^{i}(\wh C)\simeq\QQ^{i}(C)$ by construction, so~\eqref{Eq:MidRel1} and~\eqref{Eq:MidRel2} follow from Definition~\ref{Def:HodgeTwist} for $\mu$.
If $i=k$, then we have $\wh\mu(c_j) \subset \Lambda_{0,1}\perp_{\wh\Or} C^{n-i}_\rmdeg$ using \eqref{Eq:DefOfOr}, so we can assume that $c_{1} = q_\alpha$, $c_2 = q_\beta$ for some $\alpha,\beta\in I$.

In (a)--(c) we have using \eqref{Eq:DefOfMu} and \eqref{Eq:DefOfOr},
\begin{equation*}
	\la\wh\mu(q_\alpha),q_\beta\ra_{\wh\Or} = - \la\Dd x_\alpha,q_\beta\ra_{\wh\Or} = -\la q_\alpha,q_\beta\ra.
\end{equation*}
Therefore, \eqref{Eq:MidRel1} holds.

Suppose that $n$ is even and $k=n/2$ (case (d)).
If $\Char(\K)\neq 2$, then $\wh\mu$ is given by~\eqref{Eq:DefOfMuInTheMiddleDegree}, and we have
\[
	\la q_\alpha,\Dd x_\beta\ra_{\wh\Or} = \pm\la \Dd x_\beta,q_\alpha\ra_{\wh\Or} = \pm \la q_\beta,q_\alpha\ra = \la q_\alpha,q_\beta\ra,
\]
where $\pm$ is the parity of $n/2$.
Altogether we have
\begin{align*}
	\la \wh\mu(q_\alpha),q_\beta\ra_{\wh\Or} + \la q_\alpha,\wh\mu(q_\beta)\ra_{\wh\Or} = -\frac{1}{2}\bigl(\la\Dd x_\alpha,q_\beta\ra_{\wh\Or} + \la q_\alpha,\Dd x_\beta\ra_{\wh\Or}\bigr) = -\la q_\alpha,q_\beta\ra.
\end{align*}
Therefore, \eqref{Eq:MidRel2} holds.
If $\Char(\K)=2$, then $q_\alpha$ ($\alpha\in I = \{1,\dotsc,2l\}$) is a symplectic basis, and~$\wh\mu$ is given by~\eqref{Eq:DefOfMuSympl}.
Using $\la \Dd x_\alpha, q_\beta\ra = \la q_\alpha, \Dd x_\beta \ra = \la q_\alpha,q_\beta\ra$ we have for all $\alpha, \beta\in \{1,\dotsc 2l\}$ the following:
\begin{align*}
	\la \wh\mu(q_\alpha),q_\beta\ra + \la q_\alpha, \wh\mu(q_\beta)\ra &= \left.\begin{cases}
			\la \wh\mu(q_\alpha),q_{\beta}\ra  & \text{if }\alpha\in\{1,\dotsc,l\}, \beta=\alpha+l, \\
			\la q_{\alpha},\wh\mu(q_\beta)\ra  & \text{if }\beta\in\{1,\dotsc,l\}, \alpha=\beta+l, \\
		0	 & \text{otherwise},
\end{cases}\right\} \\
&= -\la q_\alpha,q_\beta\ra.
\end{align*}
Therefore, \eqref{Eq:MidRel2} holds.
\qedhere~\textit{Sublemma}
\end{proof}
Definition~\ref{Def:HodgeTwist} requires~\eqref{Eq:MidRel1} and~\eqref{Eq:MidRel2} to hold for all $c_1^\prime\in \wh{C}^i$, $c_2^\prime\in\wh C^{n-i}$, whereas the Sublemma only shows them for~$c_1\in C^i$, $c_2\in C^{n-i}$.
In order to deal with this, consider the following complement of $C$ in $\wh C$:
\begin{equation}\label{eq:complementZ}
	Z \coloneqq \wh{\pi}_\HH^\perp(\Lambda_{1,0}\cdot V)\oplus \bigoplus_{i=2}^\infty \wh C^i.
\end{equation}
So, we can write $c_1^\prime = c_1 + \xi_1$, $c_2^\prime = c_2 + \xi_2$ for unique $c_1\in C^i$, $c_2\in C^{n-i}$ and $\xi_1\in Z^i$, $\xi_2\in Z^{n-i}$.
Looking at Table~\ref{tab:coexact}, we have $\xi_1 = \xi_2 = 0$ except for the following cases:
\begin{enumerate}
	\item[(i)] $\xi_1\in \Lambda_{1,0}$ and $\xi_2=0$ for $i=k-1$ in (a),
	\item[(ii)] $\xi_1, \xi_2\in \Lambda_{1,0}$ for $i=k-1$ in (b),
	\item[(iii)] $\xi_1=0$ and $\xi_2\in \Lambda_{1,0}$ for $i=k$ in~(c).
\end{enumerate}
We have $\xi_2\in \Lambda_{1,0}\subset \wh C^{n-k}_\rmdeg\perp_{\wh \Or}\wh\mu(\wh C^k)=\Lambda_{0,1}$ in (iii) and $\xi_1\in \Lambda_{1,0} \subset \wh C^{k-1}_{\rmdeg}\subset \ker\wh\mu$ in~(i) by the definitions of $\wh\Or$ and $\wh\mu$.
Therefore, \eqref{Eq:MidRel1} holds.
In contrast to this, $\xi_1, \xi_2\in \Lambda_{1,0}$ in (ii) are not a priori perpendicular to $\wh\mu(\wh C^{n/2}) = \mu(C^{n/2})$.
However, the assumption~\eqref{Eq:AdditionalAssumption} makes~\eqref{Eq:MidRel2} trivial.
This shows that $\wh\mu\colon \QQ(\wh C) \to \im\wh\Dd$ is a Hodge twist of degree $k$.
The lemma is thereby proven under the assumption that $\ALambda = \Lambda$.

Finally, suppose that $k$ is odd and $\Char(\K)=p\neq 0$ ($\Lambda\subsetneq\ALambda$).
Consider the ideal $\Lambda^\prime\coloneqq (x_{\alpha,j},y_{\alpha,j}\mid \alpha\in I, j\in\N)\subset \ALambda$, and define $V^\prime \coloneqq \Lambda^\prime\otimes V$ and $\wt V \coloneqq \Lambda \otimes V$.
We then have $\ALambda \simeq \Lambda \oplus \Lambda^\prime$ and $\wh V \simeq \wt V\oplus V^\prime$ as vector spaces.
We define $\wh \Or$ on $\wt V$ as above and extend it to $\wh\Or\colon \wh V \to \K$ by $\wh\Or|_{V^\prime} \coloneqq 0$.
Because $\wt V_{\le 2} \subset \ker\wh\Or$ and $\Dd V^\prime\subset \wt V_{\ge 2} \oplus V^\prime$, the equation $\wh\Or \circ \wh \Dd = 0$ holds.
This together with $\wh\Or\circ\iota = \wt \Or \circ \iota = \Or$ implies that~$\wh \Or$ is the required orientation.
The requirements on the pre-Hodge decomposition $\wh V = \HH\oplus\im\wh\Dd\oplus\wh C$ and the Hodge twist $\wh\mu\colon \QQ(\wh C)\to \im\wh\Dd$ involve only degrees $i\in [n-k,k]$. 
Hence, for $i\not\in[n-k,k]$, we can pick an arbitrary complement~$\wh C^i$ of $(\ker\wh\Dd)^i$ such that $\wh C^i \perp_{\wh\Or} \HH^{n-i}$ and define $\wh\mu\coloneqq 0$ on $\QQ^{i}(\wh C)$.
As for the middle degrees, we observe that $k\ge 3$ implies
\[
		\deg x_{\alpha,1} = p(k-1) - 1 \ge 2k-3\ge k.
\]
For $k>3$ we thus have $\wh V^i = \wt V^i$ for all $i\in [0,k]$, so $\wh C^i$ and $\wh\mu\colon \QQ^{i}(\wh C)\to (\im\wh\Dd)^i$ for $i\in [n-k,k]$ can be constructed as above, and all previous arguments remain valid.
The remaining case is $p=2$ and $k=3$, when $\wh V^3 = \wt V^3\oplus\Span_\K\{x_{\alpha,1}\mid \alpha\in I\}$.
Because $V\into \wt V$ is a quasi-isomorphism in degree $3$, we can still construct the pre-Hodge decomposition $\wt V^3 = \HH^3\oplus (\im \wt \Dd)^3\oplus \wt C^3$ as above.
We have $x_{\alpha,1}\in \wh V_\rmdeg$ (since $V^\prime\subset\ker\wh\Or$), so setting $\wh C^3 \coloneqq \wt C^3\oplus\Span_\K\{x_{\alpha,1}\mid \alpha\in I\}$ gives a pre-Hodge decomposition $\wh V^3 = \HH^3\oplus(\im\wh\Dd)^3\oplus \wh C^3$ such that $\QQ^{3}(\wh C)\simeq \QQ^{3}(\wt C)$.
The Hodge twist $\wh\mu\colon \QQ^{3}(\wh C)\to (\im \wh\Dd)^3$ can thus be defined as above again.
\end{proof}

\begin{remark}\label{Rem:Consequence}
	In fact, we will not need the extension property~\eqref{eq:extension}, but rather the following conclusion of Lemmas~\ref{Lem:Extension} and~\ref{Lem:HodgeTwist}:
	\begin{equation}\label{Eq:Conclusion}
		\begin{gathered}	
			\text{If }V\text{ is of Hodge type in degrees }i\in [n-(k-1),k-1],\\
			\text{then }\wh V\text{ is of Hodge type in degrees }i\in [n-k,k].
		\end{gathered}
		\tag{{$\dagger\dagger$}}
	\end{equation}
	Note that~\eqref{Eq:AdditionalAssumption} for $k=n/2+1$ can always be satisfied by modifying $V=\HH\oplus\im\Dd\oplus C$ into a Hodge decomposition in degree $n/2$ using the given Hodge twist $\mu$.
\end{remark}

We will use the following notion of a quasi-isomorphic extension:

\begin{definition}\label{Def:Extension}
	Given an injective PDGA morphism $\iota\colon V\into V^\prime$, a \emph{PDGA retraction} $\pi\colon V^\prime \to V$ is a PDGA morphism such that $\pi \circ \iota = \Id$.
	If such $\pi$ exists, we say that $V$ is a \emph{PDGA retract} of $V^\prime$.
	We call a PDGA retract \emph{oriented} if both~$V$ and~$V^\prime$ are equipped with a chain-level orientation that is preserved by $\iota$. 
\end{definition}
\begin{remark}
Any PDGA retraction $\pi\colon V^\prime \to V$ is a PDGA quasi-isomorphism---the surjectivity of $\pi_*\colon \Hh(V^\prime)\to\Hh(V)$ follows from $\pi\circ \iota = \Id$, and the injectivity of $\pi_*$ follows from the fact that it preserves the nondegenerate pairing on $\Hh(V^\prime)$.
\end{remark}

The following extension theorem is crucial in our proof of the existence statement for 1-connected Poincar\'e duality models: 

\begin{proposition}\label{Prop:ExtensionSimplyConnected}
	Let $(V,\Dd,\Or)$ be a \emph{1-connected oriented PDGA} of degree $n\in\N_0$ over a field $\K$.
If $\Char(\K)=2$ and $n\ge 4$ is even, we additionally assume that~$V^{n/2}$ admits a Hodge decomposition.
Then there exists a connected oriented PDGA $(\wh V,\wh \Dd,\wh \Or)\supset (V,\Dd,\Or)$ of degree $n$ of Hodge type and a PDGA retraction $\pi\colon \wh V\onto V$.
Moreover, $\wh V$ is $1$-connected for $n\neq 4$, and $\wh V$ is of finite type provided that~$V$ is.
\end{proposition}

\begin{proof}
	Let $V=\HH\oplus \im \Dd\oplus C$ be a pre-Hodge decomposition of $(V,\Dd,\Or)$.
	The fact that $V$ is 1-connected implies $V^0 = \HH^0 = \Span_\K\{1\}$, $V^1=0$, and $V^2 = \HH^2\oplus C^2$.
	Therefore, $V=\HH\oplus\im\Dd\oplus C$ is a Hodge decomposition and the proposition holds trivially for $n\in [0,3]$ (the case $n=1$ is vacuous due to Poincar\'e duality on $\Hh(V)$).

	For $n=4$, it may be necessary to add generators in degree $1$ to handle $C^2\not\perp C^2$.
	The construction in the proof of Lemma~\ref{Lem:Extension} for $n=4$ and $k=2$ gives an oriented PDGA $(\wh V=\Lambda\otimes V,\wh \Dd,\wh \Or)\supset (V,\Dd,\Or)$ and a pre-Hodge decomposition $\wh V = \HH\oplus\im\wh\Dd\oplus\wh C$ (even though $n<5$). 
	We have $(\Lambda_{1,0}\cdot V)^2 = \Lambda_{1,0}\cdot V^1 = 0$ and $\wh C_2^2 = \Lambda_{2,0}\subset \wh V_\rmdeg$, so
	\[
		\wh C^2 = (\wh C_{\le 1}\oplus \wh C_2)^2 = C^2 \oplus \Lambda_{2,0}\quad\text{and}\quad\wh C^2_\rmdeg = C^2_\rmdeg \oplus \Lambda_{2,0}.
	\]
	This differs from Table~\ref{tab:coexact} (d).
	Nevertheless, we see that the inclusion $C^2\into \wh{C}^2$ induces an isomorphism $\QQ^{2}(C)\simeq \QQ^{2}(\wh C)$, so $\wh\mu\colon \QQ(\wh C)\to \im\wh\Dd$ can still be defined on $\QQ^2(\wh C)$ by~\eqref{Eq:DefOfMuInTheMiddleDegree}, which assures that~\eqref{Eq:MidRel2} with $i=k=2$ holds.
	Moreover, the complement $Z$ from~\eqref{eq:complementZ} satisfies $Z^2=\Lambda_{2,0}\subset \wh V_\rmdeg$, so~\eqref{Eq:MidRel2} still implies Definition~\ref{Def:HodgeTwist} for $\wh\mu$.
	Further, we set
	\[
		\wt C^0 \coloneqq 0,\quad\wt C^1 \coloneqq \Lambda_{1,0},\quad\wt C^2 \coloneqq \mathrm{graph}(\wh\mu)^2,\quad\text{and}\quad\wt C^i \coloneqq \wh C^i\quad\text{for all }i\ge 3.
	\]
	We see that $\wh V = \HH\oplus\im\wh\Dd\oplus\wt C$ is a Hodge decomposition in degrees $i\not\in\{1,3\}$.
	Moreover, writing $\wt C^3 = C^3 \oplus \Lambda_{1,0}\cdot V^2\oplus \Lambda_{1,1}\oplus\Lambda_{3,0}$, we see that $\wt C^1\perp_{\wh \Or} \wt C^3$, so $\wh V = \HH\oplus\im\wh\Dd\oplus\wt C$ is a Hodge decomposition in degrees $i\in\{1,3\}$ as well.

	For $n\ge 5$ we apply Lemma~\ref{Lem:Extension} with $k=\lceil n/2\rceil,\dotsc,n-2$ iteratively, starting with $(V,\Dd,\Or)$ together with a pre-Hodge decomposition $V=\HH\oplus\im\Dd\oplus C$, which exists by Lemma~\ref{Lem:PreHodge}.
	After the iterations, we obtain an oriented PDGA $(\wh V,\wh \Dd,\wh \Or)\supset (V,\Dd,\Or)$, which is of Hodge type in degrees $i\in[2,n-2]$ by~\eqref{Eq:Conclusion}, and in degrees $i\in \{0,1,n-1,n\}$ due to 1-connectedness.
	Therefore, $(\wh V,\wh\Dd,\wh \Or)$ is of Hodge type.

	In all cases, the DGA $\wh V$ is isomorphic to the tensor product of $V$ with a free acyclic CDGA $\ALambda$ (constructed iteratively).
	Let $V^\prime\coloneqq V\otimes \ALambda^{>0}$, so that $\wh V \simeq V \oplus V^\prime$ as vector spaces.
	Because~$V$ is a dg-subalgebra and $V^\prime$ a dg-ideal in $\wh V$, the projection $\pi\colon \wh V \onto V$ is a DGA morphism.
	Because $\iota\colon V\into \wh V$ is an orientation preserving DGA quasi-isomorphism and we have $\pi\circ\iota=\Id$ by construction, $\pi$ is a PDGA retraction according to Definition~\ref{Def:Extension}.
\end{proof}

We refer to an oriented PDGA $\wh V$ of Hodge type such that $V$ is an oriented PDGA retract of $\wh V$ as an \emph{extension of Hodge type} of $V$.
In Example~\ref{Ex:NoHodgeExt} we construct for any field $\K$ with $\Char(\K)=2$ and any even $n\ge 4$ a 1-connected dPD algebra of degree~$n$ over~$\K$ that does not admit an extension of Hodge type.
In Example~\ref{Ex:Degree4} we construct for any field~$\K$ with $\Char(\K)\neq 2$ a 1-connected PDGA of degree $n=4$ over $\K$ that does not admit a \emph{1-connected} extension of Hodge type.

The next result is a version of Proposition~\ref{Prop:ExtensionSimplyConnected} for connected PDGAs.
Owing to this result, we can drop the assumption $\Hh^2=0$ in the uniqueness statement for 1-connected Poincar\'e duality models if $n$ is odd.
Unlike Proposition~\ref{Prop:ExtensionSimplyConnected}, we do not construct a Hodge decomposition in the middle degree $n/2$ if~$n$ is even, thus we make no assumptions on $\K$.

\begin{proposition}\label{Prop:ExtensionConnected}
	Let $(V,\Dd,\Or)$ be a \emph{connected} oriented PDGA of degree $n\in\N_0$ \emph{of finite type} over a field $\K$.
	For even $n\in \N$ we additionally assume that~$V^{n/2}$ admits a Hodge decomposition.
	Then there is an oriented PDGA $(\wh V,\wh\Dd,\wh\Or)\supset (V,\Dd,\Or)$ of Hodge type of finite type that retracts onto $V$ as a PDGA.
\end{proposition}

\begin{proof}
	The proposition holds trivially for $n\in\{0,1\}$ for degree reason.
	The case $n=2$ is implied by the additional assumption. 

	Suppose that $n\ge 3$.
	We will prove a modification of Lemma~\ref{Lem:Extension} (\emph{the Lemma}) that allows for $V^1\neq 0$ and $n\ge 3$, but that additionally assumes that~$V$ is of finite type and $k>n/2$.
	We are only interested in the property~\eqref{Eq:Conclusion} of~$\wh V$, so we can ignore the assumption~\eqref{Eq:AdditionalAssumption} for even $n$ and $k=n/2+1$ since it can always be satisfied by the choice of a pre-Hodge decomposition.
	Because $n\ge 3$, the least odd $k>n/2$ is $k=3$, so we can suppose that $\Char(\K)=0$ or $k\in 2\N$, so that $\ALambda = \Lambda$ and $\wh V = \Lambda \otimes V$.
	The case when $\Char(\K)=p\neq 0$ and $k$ is odd can be dealt with as in the last paragraph of the proof of the Lemma.

	An inspection of the proof of the Lemma shows that the construction of the oriented PDGA $(\wh V,\wh \Dd,\wh \Or)\supset (V,\Dd,\Or)$ that PDGA retracts onto $(V,\Dd,\Or)$ and the construction of the pre-Hodge decomposition $\wh V=\HH\oplus\im\wh\Dd\oplus\wh C$ are valid for all $n\in \N_0$ regardless of $V^1\neq 0$.	
	Because $k>n/2$, we are in one of the cases (a)--(c) from Table~\ref{tab:coexact}.
	However, the table has to be corrected as follows:
	\begin{enumerate}
		\item[(i)] by adding to~$\wh{C}^k$ the subspace $\wh\pi_{\HH}^\perp(\Lambda_{1,0}\cdot V^1)$ (since $V^1\neq 0$),
		\item[(ii)] by adding to $\wh{C}^k$ the subspace~$\Lambda_{2,0}$ for $(n,k)=(3,2)$ (since $2(k-1)=k$),
		\item[(iii)] by modifying~$\wh{C}^k_\rmdeg$ and~$\wh{C}^{n-k}_\rmdeg$ accordingly.
	\end{enumerate}
	On the other hand, the components $\wh{C}^i$ for $i<k$ remain unchanged.
	We conclude that we can still construct a Hodge twist $\wh\mu\colon\QQ(\wh C)\to \im\wh\Dd$ of degree $k-1$ such that~\eqref{eq:extension} holds (relevant in (a) and (b)); however, we may not be able to extend $\wh\mu$ to a Hodge twist of degree $k$ as in the Lemma because we may have $\QQ^k(\wh C)\not\simeq\QQ^k(C)$ due to (i)--(iii) (relevant in (a)--(c)).	

	Our strategy to over come the above problem is to reapply the construction for \emph{the same~$k$} iteratively.
	We start with $(V^{(0)},\Dd^{(0)},\Or^{(0)})=(V,\Dd,\Or)$ and $V=\HH\oplus\im\Dd\oplus C = \HH\oplus \im\Dd^{(0)}\oplus C^{(0)}$, and we obtain for each $j\in\N$ an oriented PDGA 
	\[ (V^{(j)}=\Lambda^{(j)}\otimes V^{(j-1)}, \Dd^{(j)},\Or^{(j)}) \]
	with a pre-Hodge decomposition $V^{(j)}=\HH\oplus\im\Dd^{(j)}\oplus C^{(j)}$ and a Hodge twist $\mu^{(j)}\colon \QQ(C^{(j)})\to \im\Dd^{(j)}$ of degree $k-1$ that satisfies~\eqref{eq:extension} with respect to $\mu^{(j-1)}$.
	We will argue that there is some $j_0\in\N$ such that
	\begin{equation}\label{Eq:EndCondition}
		\QQ^{k}(C^{(j_0)})\simeq \QQ^{k}(C^{(j_0-1)})
	\end{equation}
	and will proceed with $V^{(j_0)}$ as in the 1-connected case.
	We denote by $\la-,-\ra^{(j)}$ the pairing on $V^{(j)}$ induced by $\Or^{(j)}$, by $\perp^{(j)}$ the corresponding perpendicularity relation, and by $\pi_{\HH}^{(j)}\colon V^{(j)}\onto \HH\subset V^{(j)}$ the harmonic projection with respect to $\la-,-\ra^{(j)}$.

	Table~\ref{tab:coexact2} shows the corrected components $C^{(j)n-k}$ and $C^{(j)k}$ obtained from~\eqref{eq:defofcoexact} and~\eqref{eq:defofcoexact2} by applying the modifications (i)--(iii).
	We now have four cases (a), (b), (c$^\prime$), and (c$^\dprime$).%
	\begin{table}
	\caption{The coexact part in degrees $k$ and $n-k$ after the $j$-th iteration.}\label{tab:coexact2}
	\renewcommand{\descriptionlabel}[1]{\hspace{\labelsep}\normalfont\underline{#1}:}
	\begin{description}[itemindent=-.35cm]
		\item[(a), (b)]
	\begin{align*}
		C^{(j) n-k} &= C^{n-k}\\
		C^{(j) k} &= C^{k} \oplus \pi_{\HH}^{\perp(1)}(\Lambda_{1,0}^{(1)}\cdot V^1)\oplus \dotsb \oplus \pi_{\HH}^{\perp(j)}(\Lambda^{(j)}_{1,0}\cdot V^1)
	\end{align*}
	\item[(c$^\prime$) $\equiv$ (c) with $(n,k)\neq (3,2)$]
		\begin{align*}
		 C^{(j) k - 1} &= C^{k-1} \oplus \Lambda^{(1)}_{1,0}\oplus \dotsb \oplus \Lambda^{(j)}_{1,0} \\
		 C^{(j) k} &= C^{k} \oplus \pi_{\HH}^{\perp(1)}(\Lambda_{1,0}^{(1)}\cdot V^{1})\oplus \dotsb \oplus \pi_{\HH}^{\perp(j)}(\Lambda^{(j)}_{1,0}\cdot V^{1})
		\end{align*}
	\item[(c$^{\dprime}$) $\equiv$  (c) with $(n,k)=(3,2)$]
	\begin{align*}
		C^{(j)1} &= C^1 \oplus \Lambda^{(1)}_{1,0}\oplus \dotsb \oplus \Lambda^{(j)}_{1,0} \\
		C^{(j)2} &= C^2 \oplus \bigoplus_{a= 1}^j \pi_{\HH}^{\perp(a)}(\Lambda_{1,0}^{(a)}\cdot V^{(a-1)1})\oplus\bigoplus_{a=1}^j \Lambda_{2,0}^{(a)}
	\end{align*}
	\end{description}
	\end{table}
	Using that the restriction of $\la-,-\ra^{(j)}$ to $V^{(j-1)}\subset V^{(j)}$ equals $\la-,-\ra^{(j-1)}$, that $C^{(j-1)}\subset C^{(j)}$, that $C^{(j)n-k} = C^{(j-1)n-k}$ in (a) and (b), and that $C^{(j)n-k} = C^{(j-1)n-k}\oplus \Lambda^{(j)}_{1,0}$ with $\Lambda^{(j)}_{1,0}\perp^{(j)} C^{(j-1)k}$ in (c$^\prime$) and (c$^{\dprime}$), we obtain
	\begin{equation}\label{eq:perpcond}
		C^{(j-1)k}_\rmdeg = C^{(j-1)k}\cap C^{(j)k}_\rmdeg.
	\end{equation}
	This implies that the inclusion $C^{(j-1)k}\into C^{(j)k}$ induces an injection $\QQ^{k}(C^{(j-1)})\into\QQ^{k}(C^{(j)})$, so the sequence $\dim \QQ^{k}(C^{(j)})$ ($j\in \N_0$) is monotone increasing. 
	The following Sublemma implies that the sequence $\dim \QQ^{n-k}(C^{(j)})$ ($j\in \N_0$) is bounded. 
	Since~$V$ is of finite type, $C^{(j)}$ is also of finite type, so we have Poincar\'e duality $\QQ^k(C^{(j)})\simeq \QQ^{n-k}(C^{(j)})$.
	Therefore, the two sequences agree, and the monotonicity and boundedness imply the existence of $j_0\in\N$ such that~\eqref{Eq:EndCondition} holds.

	\begin{sublemma*}
	The composition $C^{n-k}\into C^{(j)n-k}\onto\QQ^{n-k}(C^{(j)})$ is surjective.	
	\end{sublemma*}

	\begin{proof}
	In (a) and (b) we have $C^{(j)n-k}=C^{n-k}$, so the statement holds trivially.
	Looking at Table~\ref{tab:coexact2}, the statement in (c$^\prime$) and (c$^{\dprime}$) will follow if we show that
	\begin{equation}\label{eq:surjectivity}
		\Lambda_{1,0}^{(1)}\oplus\dotsb\oplus\Lambda_{1,0}^{(j)}\subset C^{(j)k-1}_\rmdeg.
	\end{equation}

	In (c$^\prime$), the definition of $\Or^{(a)}$ gives for each $a\in [1,j]$,
	\[
		\Lambda_{1,0}^{(a)}\perp^{(a)}\HH^{(a)k}\oplus C^{(a)k} = \HH^k\oplus C^k \oplus\bigoplus_{b=1}^a \pi_{\HH}^{\perp(b)}(\Lambda_{1,0}^{(b)}\cdot V^1),
	\]
	so we have $\Lambda_{1,0}^{(a)}\perp^{(j)} \HH^k\oplus C^k$ and $\Lambda_{1,0}^{(a)}\perp^{(j)}\pi_\HH^{\perp(b)}(\Lambda_{1,0}^{(b)}\cdot V^1)$ for all $b\in[1,a]$.
	For each $b\in(a,j]$ the cyclicity of $\la-,-\ra^{(j)}$ with respect to $\cdot$ (Property \eqref{Eq:Product}) gives
	\begin{equation}\label{eq:equivalences1}
	\begin{aligned}
		\Lambda_{1,0}^{(a)}\perp^{(j)} \pi_\HH^{\perp(b)}(\Lambda_{1,0}^{(b)}\cdot V^{1}) 
		\quad &\Longleftrightarrow\quad\Lambda_{1,0}^{(a)}\perp^{(j)} \Lambda^{(b)}_{1,0}\cdot V^1 \\
		&\Longleftrightarrow\quad\Lambda_{1,0}^{(b)}\perp^{(j)} \Lambda_{1,0}^{(a)}\cdot V^1 \\
		&\Longleftrightarrow\quad\Lambda_{1,0}^{(b)}\perp^{(j)}\pi_\HH^{\perp(a)}(\Lambda_{1,0}^{(a)}\cdot V^1).
	\end{aligned}
	\end{equation}
	We see that~\eqref{eq:surjectivity} holds.

	In (c$^\dprime$), the definition of $\Or^{(a)}$ gives for each $a\in[1,j]$ the following: $\Lambda_{1,0}^{(a)}\perp^{(j)}\HH^2\oplus C^2$, $\Lambda_{1,0}^{(a)}\perp^{(j)}\pi_\HH^{\perp(b)}(\Lambda_{1,0}^{(b)}\cdot V^{(b-1)1})$, and $\Lambda_{1,0}^{(a)}\perp^{(j)}\Lambda_{2,0}^{(b)}$ for all $b\in[1,a]$.
	Writing 
	\[
		V^{(b-1)1} = V^1\oplus \bigoplus_{c=1}^{b-1} \Lambda_{1,0}^{(c)},
	\]
	we conclude that $\Lambda_{1,0}^{(a)}\perp^{(j)}\pi_\HH^{\perp(b)}(\Lambda_{1,0}^{(b)}\cdot V^{(b-1)1})$ is equivalent to $\Lambda_{1,0}^{(a)}\perp^{(j)} \Lambda_{1,0}^{(b)}\cdot V^1$ and $\Lambda_{1,0}^{(a)}\perp^{(j)} \Lambda_{1,0}^{(b)}\cdot\Lambda_{1,0}^{(c)}$ for all $c\in[1,b)$.
	Suppose now that $b\in (a,j]$.
	The definition of~$\Or^{(b)}$ gives $\Lambda_{2,0}^{(b)}\perp^{(j)}\Lambda_{1,0}^{(a)}$.
	By~\eqref{eq:equivalences1} we have $\Lambda_{1,0}^{(a)}\perp^{(j)} \Lambda_{1,0}^{(b)}\cdot V^1\Longleftrightarrow\Lambda_{1,0}^{(b)}\perp^{(j)} \Lambda_{1,0}^{(a)}\cdot V^1$.
	The cyclicity of $\la-,-\ra^{(j)}$ with respect to $\cdot$ gives
	\begin{align*}
		\Lambda_{1,0}^{(a)}\perp^{(j)}\Lambda_{1,0}^{(b)}\cdot\Lambda_{1,0}^{(c)}
		\quad\Longleftrightarrow\quad\Lambda_{1,0}^{(b)} \perp^{(j)}\Lambda_{1,0}^{(c)}\cdot\Lambda_{1,0}^{(a)} 
		\quad\Longleftrightarrow\quad\Lambda_{1,0}^{(b)} \perp^{(j)}\Lambda_{1,0}^{(a)}\cdot\Lambda_{1,0}^{(c)}.
	\end{align*}
	The middle perpendicularity relation was already shown for $1\le a<c \le b$, the right one was shown for $1\le c < a \le b$, and the relation for $a=c$ is equivalent to $\Lambda_{1,0}^{(b)}\perp^{(j)}\Lambda_{2,0}^{(a)}$, which was also shown.
	We see that $\Lambda_{1,0}^{(a)}\perp^{(j)}\Lambda_{1,0}^{(b)}\cdot\Lambda_{1,0}^{(c)}$ holds also for all $b\in(a,j]$ and $c\in[1,b)$.
	We conclude that~\eqref{eq:surjectivity} holds.
	\qedhere~\textit{Sublemma}
	\end{proof}	

	Having a $j_0\in \N$ such that~\eqref{Eq:EndCondition} holds (we may assume that $j_0$ is the minimal natural number with this property), we can define $\mu^{(j_0)}$ on $\QQ^{k}(C^{(j_0)})\simeq\QQ^{k}(C^{(j_0-1)})$ by~\eqref{Eq:DefOfMu} (with $q_\alpha^{(j_0)}$ and $x_{\alpha}^{(j_0)}$ for $\alpha\in I^{(j_0)}$) to ensure that~\eqref{Eq:MidRel1} holds for all $c_1\in C^{(j_0-1)k}$ and $c_2\in C^{(j_0-1)n-k}$ (with respect to $\la-,-\ra^{(j_0)}$).
	In order to deduce that~$\mu^{(j_0)}$ satisfies Definition~\ref{Def:HodgeTwist} of a Hodge twist of degree $k$, we have to show that~\eqref{Eq:MidRel1} holds also for $c_1^\prime = c_1 + \xi_1$ and $c_2^\prime = c_2 + \xi_2$, where $\xi_1\in Z^{(j_0)k}$ and $\xi_2\in Z^{(j_0)n-k}$ are elements of the complement $Z^{(j_0)}$ of $C^{(j_0-1)}$ in $C^{(j_0)}$ defined in~\eqref{eq:complementZ}.
	We distinguish the following cases:
	\begin{enumerate}
		\item[(i)] $\xi_1\in \pi_{\HH}^{\perp(j_0)}(\Lambda_{1,0}^{(j_0)}\cdot V^1)$ and $\xi_2=0$ in (a) and (b),
		\item[(ii)] $\xi_1\in \pi_{\HH}^{\perp(j_0)}(\Lambda_{1,0}^{(j_0)}\cdot V^1)$ and $\xi_2\in \Lambda_{1,0}^{(j_0)}$ in (c$^\prime$),
		\item[(iii)] $\xi_1\in \pi_{\HH}^{\perp(j_0)}(\Lambda_{1,0}^{(j_0)}\cdot V^{(j_0-1)1})$ and $\xi_2\in \Lambda_{1,0}^{(j_0)}$ in (c$^\dprime$).
	\end{enumerate}
	In each (i)--(iii) we have $\xi_1\in C^{(j_0)}_\rmdeg$ by~\eqref{Eq:EndCondition} and $\xi_2\in \Lambda_{1,0}^{(j_0)}\subset C^{(j_0)k-1}_\rmdeg$, so Definition~\ref{Def:HodgeTwist} indeed holds.
	Having the Hodge twist $\wh\mu$ of degree $k$, Lemma~\ref{Lem:HodgeTwist} implies that~$\wh V$ is of Hodge type in degrees $i\in [n-k,k]$.
	Therefore, \eqref{Eq:Conclusion} holds, and the modified Lemma is thereby proven.

	Having the modified Lemma, we can apply it iteratively with $k=\lceil n/2\rceil, \dotsc, n-1$ for odd $n$ and with $k=n/2+1, \dotsc, n-1$ for even $n$ and argue as in the proof of Proposition~\ref{Prop:ExtensionSimplyConnected} for $n\ge 5$.
	This time, $\wh V$ is generally only connected, so we may have to iterate until $k=n-1$.
	Moreover, we also cover the cases $(n,k)\in \{(3,2), (4,3)\}$.
	Other arguments, including the existence of a PDGA retraction, are the same. 
\end{proof}

\section{Differential Poincar\'e duality models}\label{Sec:5}

We consider the following equivalence relation on PDGAs: 
\begin{definition}
We say that two PDGAs $(V,\Dd,\Or_*)$ and $(V^\prime,\Dd^\prime,\Or_*^\prime)$ are \emph{weakly homotopy equivalent as PDGAs} if there is a zig-zag 
\begin{equation}\label{Eq:ZigZag}
	\begin{tikzcd}[column sep=2em]
	& V_1 \ar[dl,"f_0"'] \ar[dr,"f_1"] & & V_3 \cdots \ar[dl,"f_2"'] & & V_k \ar[dl,"f_{k-1}"'] \ar[dr,"f_{k}"] \\
	V & & V_2 & & \cdots V_{k-1} & & V^\prime,
\end{tikzcd}
\end{equation}
where $(V_i,\Dd_i,\Or_{i*})$ are PDGAs and $f_i$ PDGA quasi-isomorphisms.
\end{definition}

\begin{definition}\label{Def:dPDModel}
	A \emph{differential Poincar\'e duality model (dPD model)} of a PDGA $(V,\Dd,\Or_*)$ is a dPD algebra $(V^\prime,\Dd^\prime,\Or^\prime)$ such that $(V,\Dd,\Or_*)$ and $(V^\prime,\Dd^\prime,\Or_{*}^\prime)$ are weakly homotopy equivalent as PDGAs.
\end{definition}

In contrast to this, the definition of a dPD model in \cite{Lambrechts2007} only requires $(V,\Dd)$ and $(V^\prime,\Dd^\prime)$ to be weakly homotopy equivalent \emph{as DGAs}---that is, there is a zig-zag~\eqref{Eq:ZigZag} with DGAs $(V_i,\Dd_i)$ and DGA quasi-isomorphisms~$f_i$.
Weak homotopy equivalence of PDGAs clearly implies weak homotopy equivalence of DGAs.
We illustrate in Example~\ref{Ex:CP2} that the reverse implication does not hold; nevertheless, the notions of \emph{formality} coincide:

\begin{remark}\label{Rem:Formality}
Recall that a DGA $(V,\Dd)$ is called \emph{formal} if it is weakly homotopy equivalent (as a DGA) to its cohomology $(\Hh(V),\Dd\coloneqq 0)$.
Suppose that $(V,\Dd,\Or_*)$ is a PDGA such that the underlying DGA $(V,\Dd)$ is formal, and let~\eqref{Eq:ZigZag} be a weak homotopy equivalence of DGAs $(V,\Dd)$ and $(V^\prime\coloneqq \Hh(V),\Dd^\prime\coloneqq 0)$.
Starting with $\Or_{1*} \coloneqq \Or_*\circ f_{0*}\colon \Hh(V_1)\to\K$ for $i=1$, we define the orientation $\Or_{i*}\colon \Hh(V_i)\to\K$ for each $i=2,\dotsc,k$ inductively such that $(V_i,\Dd_i,\Or_{i*})$ becomes a PDGA and~$f_{i-1}$ becomes a PDGA quasi-isomorphism.
Finally, we replace~$f_k\colon V_k \to V^\prime = \Hh(V)$ with the composition $f_k^\prime \coloneqq f_{0*} \circ f_{1*}^{-1} \circ \dotsb \circ f_{k-1 *}\circ (f_{k*})^{-1}\circ f_k\colon V_k \to V^\prime$, so that~\eqref{Eq:ZigZag} becomes a weak homotopy equivalence of PDGAs.
We conclude that $(\Hh(V),\Or_*)$ is a dPD model of a PDGA $(V,\Dd,\Or_*)$ if and only if the underlying DGA $(V,\Dd)$ is formal.

We recall the well-known result from \cite{Formality} which states that the de Rham algebra $(\Om(X),\Dd)$ of a connected closed $n$-manifold $X$ is formal if $n\le 6$. 
\end{remark}

Suppose that $(V,\Dd,\Or)$ is an oriented PDGA of Hodge type whose cohomology $\Hh(V)$ is 1-connected. 
Given a special propagator $\Pp\colon V\to V$, Lemma~\ref{Lem:SmallSubalgebra} (b), (c) implies that the associated small subalgebra $\SS_{\Pp}\subset V$ is 1-connected, of finite type, and of Hodge type, and that the inclusion $\iota\colon \SS_\Pp\into V$ is an orientation preserving quasi-isomorphism.
Lemma~\ref{Lem:HodgeType} (a) implies that the quotient map $\pi_\QQ\colon \SS_{\Pp} \onto \QQ(\SS_\Pp)$ is an orientation preserving quasi-isomorphism.
The nondegenerate quotient $\QQ(\SS_{\Pp})$ is a dPD algebra because the induced pairing $\la-,-\ra_\QQ\colon \QQ(\SS_{\Pp})\times\QQ(\SS_{\Pp})\to\K$ is nondegenerate and~$\SS_{\Pp}$ is of finite type.
Therefore, $\QQ(\SS_{\Pp})$ is a 1-connected dPD model of $(V,\Dd,\Or_*)$ via the following zig-zag of orientation preserving PDGA quasi-isomorphisms:
\begin{equation}\label{Eq:HodgeModel}
 \begin{tikzcd}
 & \SS_{\Pp}\arrow[hook',"\iota"']{dl}\arrow[two heads,"\pi_\QQ"]{dr} & \\
 V & & \QQ(\SS_\Pp).
 \end{tikzcd}
\end{equation}
We illustrate in Example~\ref{Ex:SU6} that the dPD algebras $\QQ(\SS_\Pp)$ and $\QQ(\SS_{\Pp^\prime})$ for different special propagators~$\Pp$ and $\Pp^\prime$ may not be isomorphic.

If the oriented PDGA $(V,\Dd,\Or)$ is not of Hodge type but satisfies the assumptions of Propositions~\ref{Prop:ExtensionSimplyConnected} or \ref{Prop:ExtensionConnected}, we can first construct an oriented PDGA of Hodge type $(\wh V,\wh \Dd,\wh \Or)$ that PDGA retracts onto $(V,\Dd,\Or)$ and then apply the construction depicted in~\eqref{Eq:HodgeModel} to~$(\wh V,\wh\Dd,\wh\Or)$ instead.

If $(V,\Dd,\Or)$ is a general oriented PDGA whose cohomology $\Hh(V)$ is 1-connected, we follow~\cite{Lambrechts2007} and consider its \emph{minimal Sullivan model} $\rho\colon \Lambda U\to V$ equipped with the pullback orientation.
Because $\Lambda U$ is 1-connected, we can apply Proposition~\ref{Prop:ExtensionSimplyConnected} and continue with $\wh{\Lambda U}$ as in the previous paragraph.
Moreover, because~$\Lambda U$ and~$\wh{\Lambda U}$ are of finite type, $\QQ(\wh{\Lambda U})$ is a dPD algebra, so we do not have to descend to a small subalgebra.

We now formulate our version of \cite[Theorem~1.1]{Lambrechts2007}.

\begin{assumption}
In the following three propositions we assume that $\Char(\K)= 0$.
\end{assumption}

\begin{proposition}\label{Prop:Existence}
	Every PDGA $(V,\Dd,\Or_*)$ of degree $n\in\N_0$ with 1-connected cohomology admits a 1-connected dPD model $(V^\prime,\Dd^\prime,\Or^\prime)$ of the form 
	\begin{equation}\label{Eq:HodgeModelGeneral}
		 \begin{tikzcd}
			 & \wh{\Lambda U}\arrow["\rho\,\circ\,\pi"']{dl}\arrow[two heads,"\pi_\QQ"]{dr} & \\
			 V & & V^\prime \coloneqq \QQ(\wh{\Lambda U}),
		 \end{tikzcd}
	 \end{equation}
	 where $\rho\colon \Lambda U\to V$ is a minimal Sullivan model of $(V,\Dd)$, $\wh{\Lambda U}\supset \Lambda U$ is a 1-connected oriented PDGA of Hodge type of finite type, $\pi\colon \wh{\Lambda U}\onto\Lambda U$ is a PDGA retraction, and $\pi_\QQ\colon \wh{\Lambda U}\onto \QQ(\wh{\Lambda U})$ is the quotient map.
\end{proposition}

\begin{proof}
	Since $\Hh(V)$ is 1-connected, a \emph{minimal Sullivan model} $\rho\colon \Lambda U \to V$ can be constructed degreewise as in the paragraph preceding \cite[Proposition~12.2]{Felix2001}.
	The proposition then asserts that $U^0=U^1=0$ and $U^2 \simeq \Hh^2(V)$; moreover, we can see from the construction that $\Lambda U$ is of finite type since $\Hh(V)$ is.
	Let $\Or\colon V\to \K$ be an orientation of degree~$n$ on $(V,\Dd)$ that induces the given orientation~$\Or_*$ on cohomology.
	We equip $\Lambda U$ with the pullback orientation $\Or_{\Lambda U}\coloneqq {\Or} \circ \rho \colon \Lambda U\to \K$, which makes~$\Lambda U$ into an oriented PDGA and $\rho\colon\Lambda U\to V$ into an orientation preserving PDGA quasi-isomorphism.

	If $n\le 5$, then any pre-Hodge decomposition $\Lambda U= \HH \oplus \im\Dd \oplus C$, which exists by Lemma~\ref{Lem:PreHodge}, is Hodge for degree reasons.
	Therefore, we can take $\wh{\Lambda U} \coloneqq \Lambda U$.
	For $n\ge 6$ we apply Proposition~\ref{Prop:ExtensionSimplyConnected} to $(\Lambda U,\Dd,\Or_{\Lambda U})$ and obtain a 1-connected oriented PDGA $(\wh{\Lambda U},\wh{\Dd},\wh\Or_{\Lambda U})$ of Hodge type of finite type together with a PDGA retraction $\pi\colon\wh{\Lambda U}\onto \Lambda U$ (see Definition~\ref{Def:Extension}).
	The pairing $\la-,-\ra_{\QQ}$ is nondegenerate and $\wh{\Lambda U}$ is of finite type, so $\QQ(\wh{\Lambda U})$ is a dPD algebra.
	The quotient map $\pi_\QQ\colon \wh{\Lambda U} \onto \QQ(\wh{\Lambda U})$ is an orientation preserving PDGA quasi-isomorphism by Lemma~\ref{Lem:HodgeType}\,(a) because~$\wh{\Lambda U}$ is of Hodge type.
\end{proof}

An example of a PDGA that cannot be connected to its dPD model (any of them) by a single PDGA quasi-isomorphism is given in Example~\ref{Ex:CP7}.

Now we present our version of \cite[Theorem~7.1]{Lambrechts2007}.
Our proof is a reproduction of their proof, with the exception that we utilize our extension of Hodge type and additionally verify that the involved maps preserve cohomology orientation.

\begin{proposition}\label{Prop:Uniqueness}
Let $(V,\Dd,\Or)$ and $(V^\prime,\Dd^\prime,\Or^\prime)$ be 1-connected dPD algebras of degree $n\in\N_0$ that are weakly homotopy equivalent as PDGAs.
If $\Hh^2(V) \simeq \Hh^2(V^\prime) =0$, then there is a 1-connected dPD algebra~$(V^\dprime,\Dd^\dprime,\Or^\dprime)$ and PDGA quasi-isomorphisms
\begin{equation}\label{Eq:WeakUniqueness}
\begin{tikzcd}
& V^\dprime & \\
V \arrow[hook]{ur}{\iota}  & & V^\prime \arrow[hook',swap]{ul}{\iota^\prime},
\end{tikzcd}
\end{equation}
which are injective and preserve chain-level orientation.
\end{proposition}

\begin{proof}
	If $n=0$, then $V=\K$, so the proposition is satisfied trivially.
	The case $n=1$ is empty.
	For all $n\ge 2$ consider a zig-zag of PDGA quasi-isomorphisms~\eqref{Eq:ZigZag}. 

	If $n\in[2,4]$, then the 1-connectedness, together with Poincar\'e duality on chain level, implies that $\Dd \equiv \Dd^\prime \equiv 0$, so the dPD algebras $(V,\Dd,\Or)$ and $(V^\prime,\Dd^\prime,\Or^\prime)$ are canonically isomorphic to $(\Hh(V),\Dd\coloneqq 0, \Or_*)$ and $(\Hh(V^\prime),\Dd^\prime\coloneqq 0, \Or_*^\prime)$, respectively.
	The proposition is then satisfied by setting $V^\dprime\coloneqq \Hh(V)$, $\iota\coloneqq \Id$, and $\iota^\prime\coloneqq f_{0*}\circ f_{1*}^{-1}\circ\dotsb\circ f_{k-1*}\circ f_{k*}^{-1}\colon (\Hh(V^\prime),\Or_*^\prime)\to (\Hh(V),\Or_*)$.

	For $n\ge 5$ let $\rho\colon \Lambda U \to V$ be a \emph{minimal Sullivan model} of $(V,\Dd)$. 
	The pullback orientation $\Or_{\Lambda U} \coloneqq \Or\circ{\rho}\colon \Lambda U\to \K$ makes $\Lambda U$ into an oriented PDGA of degree~$n$ and $\rho\colon \Lambda U \to V$ into a PDGA quasi-isomorphism.
	By the \emph{lifting property} of Sullivan models \cite[Proposition~12.9]{Felix2001}, there is a DGA morphism $\rho_1\colon \Lambda U \to V_1$ such that the maps $f_0\circ\rho_1\colon \Lambda U \to V$ and $\rho\colon \Lambda U \to V$ are homotopic DGA morphisms.
	In particular, we have $f_{0*}\circ \rho_{1*} = \rho_{*}$, so we conclude that $\rho_1\colon \Lambda U \to V_1$ and $\rho_2\coloneqq f_1 \circ \rho_1\colon \Lambda U \to V_2$ are PDGA quasi-isomorphisms.
	We proceed inductively along the zig-zag to obtain a PDGA quasi-isomorphism $\rho^\prime \coloneqq f_k\circ\rho_k \colon \Lambda V\to V^\prime$.
	
	Having the PDGA quasi-isomorphisms $\rho\colon \Lambda U \to V$ and $\rho^\prime\colon \Lambda U \to V^\prime$, our initial attempt to get~\eqref{Eq:WeakUniqueness} would be to consider the tensor product of differential graded $\Lambda U$-modules $V\otimes_{\Lambda U} V^\prime$, extend it to a PDGA of Hodge type, and define~$V^\dprime$ as the nondegenerate quotient of this extension.
	However, $V\otimes_{\Lambda U} V^\prime$ is not necessary quasi-isomorphic to~$V$.%
	\footnote{The DGA $V\otimes_{\Lambda U} V^\prime$ is the \emph{pushout} of $\rho\colon\Lambda U\to V$ and $\rho^\prime\colon\Lambda U\to V^\prime$ in the category of differential graded $\Lambda U$-modules.
	However, while pushouts preserve isomorphisms, they do not necessarily preserve quasi-isomorphisms; specifically, the cohomology of a pushout diagram may not be a pushout diagram.}
	A solution from~\cite{Lambrechts2007} uses a \emph{minimal Sullivan model  $\phi\colon \Lambda U\otimes\Lambda U\otimes \Lambda W \to \Lambda U$ of the multiplication} ${\cdot}\colon \Lambda U \otimes \Lambda U\to \Lambda U$ (see \cite[Section~14]{Felix2001} for the definition of a Sullivan model of a morphism and its properties), which is in the case that $\Hh(V)$ is 1-connected and of finite type constructed in \cite[Section 15 (c), Example 1]{Felix2001}.
	Here $W = U[1]$ is the \emph{degree shift} of $U$ by $1$, which is a graded vector space defined by $U[1]^i=U^{i+1}$ for all $i\in \N_0$.

	We abbreviate $B\coloneqq V\otimes V^\prime$ and $B^\prime\coloneqq \Lambda U\otimes \Lambda U$ and consider the following diagram of DGA morphisms, where the horizontal arrows represent canonical inclusions:
	\begin{equation*}
		\begin{tikzcd}
			& & \Lambda U \\
			& B^\prime \arrow["\rho\otimes\rho^\prime"']{d}\arrow[hook]{r}\arrow[yshift=2pt,"{\cdot}"]{ru} & D^\prime\coloneqq B^\prime\otimes\Lambda W\mathrlap{\simeq B^\prime\otimes_{B^\prime}D^\prime}\arrow["\phi"']{u}\arrow["\rho\otimes\rho^\prime\otimes\Id"]{d}\\
			V\arrow[hook]{r} & B \arrow[hook]{r} & D\coloneqq B\otimes\Lambda W\mathrlap{\simeq B\otimes_{B^\prime}D^\prime.} 
		\end{tikzcd}
		\hphantom{\simeq B\otimes_{B^\prime}(B^\prime\otimes\Lambda W)}
	\end{equation*}
	Here $\phi\colon D^\prime\to \Lambda U$ is the minimal Sullivan model of $\cdot$ (so, the differential on $D^\prime$ is not induced from the tensor product $B^\prime\otimes \Lambda W$ despite the notation!), hence the triangle is commutative and $\phi$ is a quasi-isomorphism. 
	The square is a pushout diagram of differential graded $B^\prime$-modules. 
	According to \cite[Lemma~14.2]{Felix2001}, the fact that $\rho\otimes\rho^\prime$ is a quasi-isomorphism implies that $\rho\otimes\rho^\prime\otimes\Id$ is a quasi-isomorphism.%
	\footnote{This follows from the general fact that any \emph{relative Sullivan algebra} $D$ with base $B$ is semi-free as a differential graded $B$-module, so the tensor product $\otimes_{B} D$ preserves quasi-isomorphisms.}

	We equip $\Hh(D)$ with the orientation $\Or_{D*}\coloneqq \Or_{\Lambda U*}\circ \phi_* \circ (\rho\otimes\rho^\prime\otimes\Id)_*^{-1}$, so that $(D,\Dd,\Or_{D*})$ becomes a PDGA of degree $n$.
	Let $f\colon V\into D$ be the composition of inclusions along the lower side of the diagram.
	The K\"unneth theorem implies that the zig-zag from $\Lambda U$ to $V$ along the upper side of the diagram induces the map~$\rho_*$ on cohomology.
	Together with the commutativity of the diagram this yields $\Or_{D*}\circ f_* = \Or_*$, so $f$ is a PDGA quasi-isomorphism.
	We likewise obtain a PDGA quasi-isomorphism $f^\prime\colon V^\prime\into D$.

	Let $\Or_{D}\colon D\to\K$ be a chain-level orientation of degree $n$ on $D$ that induces the orientation~$\Or_{D *}$ on cohomology.
	We have $D \simeq V \otimes V^\prime \otimes \Lambda(U[1])$, where~$\Lambda U$ is of finite type and satisfies $U^0=U^1=0$ and $U^2 \simeq \Hh^2(V) = 0$.
	We see that $D$ is 1-connected and of finite type, so we can apply Proposition~\ref{Prop:ExtensionSimplyConnected} to obtain a 1-connected oriented PDGA $(\wh D,\wh \Dd,\wh\Or_{D})\supset (D,\Dd,\Or_D)$ of Hodge type of finite type that PDGA retracts onto~$D$.
	We conclude that the nondegenerate quotient $V^\dprime\coloneqq\QQ(\wh D)$ is a 1-connected dPD algebra, that the quotient map $\pi_\QQ\colon \wh D \onto \QQ(\wh D)$ is an orientation preserving PDGA quasi-isomorphism by Lemma~\ref{Lem:HodgeType}\,(a), and that the compositions $\iota\coloneqq\pi_\QQ\circ f\colon V\to V^{\dprime}$ and $\iota^\prime\coloneqq\pi_\QQ\circ f^\prime\colon V^\prime\to V^{\dprime}$ are PDGA quasi-isomorphisms.
	Moreover, $\iota$ and~$\iota^\prime$ preserve chain-level orientation and are injective by Lemma~\ref{Lem:dPDMorphism}.
\end{proof}

We show in Example~\ref{Ex:NonUniqueSimplyConnectedPDModel} that the proposition does not hold without the additional assumption $\Hh^2(V) \simeq \Hh^2(V^\prime) = 0$ if we require $V^{\dprime}$ to be 1-connected.
Nevertheless, a conjecture at the end of~\cite{Lambrechts2007} suggests that the additional assumption can be dropped if we allow $V^{\dprime 1} \neq 0$.
Using Proposition~\ref{Prop:ExtensionConnected} we can prove this for \emph{odd} $n$:

\begin{proposition}\label{Prop:LSConjectureOddCase}
	The statement of Proposition~\ref{Prop:Uniqueness} holds without the assumption $\Hh^2(V)\simeq \Hh^2(V^\prime)=0$ provided that $n$ is \emph{odd} and we allow $V^{\dprime 1}\neq 0$.
\end{proposition}

\begin{proof}
	The proof for $\Hh^2\neq 0$ goes along the same lines as the proof of Proposition~\ref{Prop:Uniqueness}.
	The only difference is in the last paragraph, where $D = V\otimes V^\prime \otimes \Lambda(U[1])$ is now only connected since $U[1]^1 = U^2 \simeq \Hh^2(V)$.
	Therefore, instead of applying Proposition~\ref{Prop:ExtensionSimplyConnected}, we need to apply Proposition~\ref{Prop:ExtensionConnected}, which requires that $n$ is odd.
	The conclusion of the proof still holds even though $\wh D$ and $V^\dprime$ are only connected.
\end{proof}

\begin{remark}\label{rem:assumption}
	In addition to odd $n$, Proposition~\ref{Prop:LSConjectureOddCase} also holds for $n\in\{0,2,4\}$ for degree reasons.
	More generally, the proposition holds for even $n\ge 6$ provided that the PDGA $D$ admits an extension of Hodge type in degree $n/2$ (cf.~Proposition~\ref{Prop:ExtensionConnected}).
\end{remark}

\begin{remark}\label{rem:assumptiononchar} 
	We could drop the assumption $\Char(\K)=0$ provided that the results on Sullivan models that we used in the proofs of Proposition~\ref{Prop:Existence} for $n\ge 6$ and Proposition~\ref{Prop:Uniqueness} for $n\ge 5$ hold.%
	\footnote{They explicitly assume that $\Char(\K)=0$ at the beginning of each \cite[Section~12--15]{Felix2001}, which is our reference for the theory of Sullivan models.
	However, upon inspection it appears that the proof of~\cite[Proposition~12.2]{Felix2001} works unmodified for any field~$\K$, and that the proof of~\cite[Proposition~12.9]{Felix2001} works for any~$\K$ if we replace the free CDGA $\Lambda$ with the free acyclic CDGA $\ALambda$ from Definition~\ref{Def:CharDef}.
	The same replacement $\Lambda\mapsto\ALambda$ would have to be done in \cite[Section 15 (c), Example 1]{Felix2001}, which would modify the expression for~$D$ in our proofs.
	However, such $D$ would still be 1-connected and of finite type, so our arguments would work.}
	However, we would still need the assumption $\Char(\K)\neq 2$ for even $n$ in order to apply Proposition~\ref{Prop:ExtensionSimplyConnected} to~$\Lambda U$ in the proof of Proposition~\ref{Prop:Existence} and to~$D$ in the proof of Proposition~\ref{Prop:Uniqueness}.
	Recall that the origin of this assumption is that an extension of Hodge type may not exist in degree~$n/2$ if $\Char(\K)=2$ (see Example~\ref{Ex:NoHodgeExt}).
	We could drop this assumption if all elements of $\QQ^{n/2}(\Lambda U)$, resp., $\QQ^{n/2}(D)$ squared to zero (cf.~\eqref{Eq:SymplecticBasis}).	
	Note that at the beginning of \cite[Section~2]{Lambrechts2007} they explicitly assume that elements of odd degree in their CDGAs square to zero even when $\Char(\K)=2$.
\end{remark}

\section{Counterexamples}

\begin{example}[Poincar\'e duality but no Hodge decomposition in characteristic 2]\label{Ex:NoHodge}
	Let $\K$ be a field of characteristic $2$.
	For \emph{any} $k\in\N$ consider the cyclic cochain complex $(V,\Dd,\la-,-\ra)$ of degree $n\coloneqq 2k$ given by $V\coloneqq\Span_\K\{1,a,\Dd a,b,\Dd b, v\}$ with $\deg a = k-1$, $\deg b =k$, $\deg v = n$, and $\la v,1\ra = \la a, \Dd b\ra=\la b,b\ra= 1$.
	Suppose that $V=\HH\oplus\im\Dd\oplus C$ is a pre-Hodge decomposition. 	
	We have $\HH^k=0$, $(\im \Dd)^k = \Span_\K\{\Dd a\}$, and there are $\alpha,\beta\in \K$ with $\beta\neq 0$ such that $C^k = \Span_\K\{ c\coloneqq \alpha \Dd a + \beta b\}$.
	We compute
	\begin{equation}\label{Eq:Characteristic2}
			\la c,c\ra = \alpha^2\la \Dd a, \Dd a\ra + \beta^2 \la b, b\ra + \alpha \beta (\la \Dd \alpha, b\ra + \la b,\Dd a\ra) = \beta^2\neq 0.\\
	\end{equation}
	We see that $C^k \not\perp C^k$, so $V=\HH\oplus\im\Dd\oplus C$ is not a Hodge decomposition.
	We conclude that $(V,\Dd,\la-,-\ra)$, despite being a cyclic cochain complex, is not of Hodge type (cf.~Lemma~\ref{Lem:CyclicIsHodge}).
\end{example}

\begin{example}[1-connected dPD algebra with no extension of Hodge type in characteristic 2]\label{Ex:NoHodgeExt}
	The cyclic cochain complex $(V,\Dd,\la-,-\ra)$ from Example~\ref{Ex:NoHodge} can be made into a dPD algebra $(V,\Dd,\Or)$ by defining a graded commutative multiplication such that $a \cdot \Dd b = b \cdot b = v$ and $\Dd a\cdot \Dd a = \Dd a\cdot b = 0$, and an orientation such that $\Or(v)=1$ (note that $(V,\Dd)$ is the quotient of the free CDGA $\Lambda_\K(a,b,\Dd a,\Dd b)$ by the dg-ideal generated by $a^2$, $ab$, $a\Dd a$, $b^2- a \Dd b$ and all elements of degree $>n$). 
	Suppose that there is an oriented PDGA $(V^\prime,\Dd^\prime,\Or^\prime)$ of degree $n$ and a PDGA quasi-isomorphism $\iota\colon V\to V^\prime$.
	Lemma~\ref{Lem:dPDMorphism} asserts that $\iota$ preserves chain-level orientation and is injective.
	Denoting $b^\prime \coloneqq \iota(b)$, we obtain $\la b^\prime,b^\prime\ra^\prime = \la b,b\ra = 1$.
	Suppose that $V^\prime = \HH^\prime\oplus \im \Dd^\prime\oplus C^\prime$ is a pre-Hodge decomposition.
	Because $\Hh^k(V^\prime)\simeq \Hh^{k}(V) = 0$, there exist $c^\prime, c^\dprime\in C^\prime$ such that $b^\prime = c^\prime + \Dd c^\dprime$.
	We compute $\la c^\prime,c^\prime \ra^\prime = 1$ as in \eqref{Eq:Characteristic2}, so $V^\prime=\HH^\prime\oplus\im\Dd^\prime\oplus C^\prime$ is not a Hodge decomposition. 
	We conclude that $(V,\Dd,\Or)$, despite being 1-connected for $n\ge 4$, does not admit an extension of Hodge type (cf.~Proposition~\ref{Prop:ExtensionSimplyConnected}).
\end{example}

\begin{example}[1-connected PDGA of degree 4 with no 1-connected extension of Hodge type] 
\label{Ex:Degree4}
Let $\K$ be a field with $\Char(\K)\neq 2$.
For \emph{even} $k\in \N$ consider the oriented PDGA $(V,\Dd,\Or)$ of degree $n\coloneqq 2k$ given by $V\coloneqq\Span_\K\{1,a,\Dd a, a^2\}$ with $\deg a = k$ and $\Or(a^2)=1$ (note that $(V,\Dd)$ is the quotient of the free CDGA $\Lambda_\K(a,\Dd a)$ by the dg-ideal generated by all elements of degree~$>n$).
Suppose that there is an oriented PDGA $(V^\prime,\Dd^\prime,\Or^\prime)$ of degree~$n$ of Hodge type and a PDGA quasi-isomorphism $\iota\colon V\to V^\prime$.
Denoting $a^\prime \coloneqq \iota(a)$, we obtain $\la a^\prime,a^\prime\ra = \la a,a\ra = 1$ as in~\eqref{eq:OrientationInTopDegree}.
Let $V^\prime = \HH^\prime\oplus\im\Dd^\prime\oplus C^\prime$ be a Hodge decomposition.
Because $\Hh^k(V^\prime)\simeq \Hh^k(V)=0$, there exist $c^\prime, c^\dprime\in C^{\prime}$ such that $a^\prime = c^\prime + \Dd c^\dprime$. 
We compute $0=\la c^\prime,c^\prime\ra^\prime = 1 - 2\la b^\prime,\Dd c^\dprime\ra^\prime$ as in~\eqref{Eq:Characteristic2}, which implies $c^\dprime \neq 0$.
Because $\deg c^\dprime = \deg a - 1 = k-1$, we get that $V^{\prime 1}\neq 0$ if $k=2$.
We conclude that $(V,\Dd,\Or)$ for $n=4$, despite being 1-connected, does not admit a \emph{1-connected} extension of Hodge type (cf.~Proposition~\ref{Prop:ExtensionSimplyConnected}).
\end{example}

\begin{example}[A PDGA that is not PDGA quasi-isomorphic to itself with the opposite orientation]\label{Ex:CP2}
Let $\K$ be a field and $n\in \N$ a multiple of $4$.
Let $\Lambda\coloneqq\Lambda_\K(a)$ be the free graded commutative $\K$-algebra with $\deg a = k\coloneqq n/2$.
Consider the quotient algebra $V\coloneqq \Lambda/(a^3) \simeq \Span_\K\{1,a,a^2\}$.
For each $\lambda\in \K\backslash\{0\}$ we have an orientation $\Or_\lambda\colon V \to \K$, $a^2\mapsto \lambda$ that makes $(V,\Or_\lambda)$ into a Poincar\'e duality algebra of degree~$n$.
For $\lambda, \lambda^\prime\in\K\backslash\{0\}$ the PD algebras $(V,\Or_\lambda)$ and $(V,\Or_{\lambda^\prime})$ are isomorphic if and only if $\lambda/\lambda^\prime$ is a square in~$\K$ (there are then two possible isomorphisms $(V,\Or_{\lambda})\to (V,\Or_{\lambda^\prime})$, $a\mapsto \pm \sqrt{\lambda/\lambda^\prime}a$).

Consider the complex projective plane $X\coloneqq \CP^2$ as a (canonically) oriented closed $4$-manifold, and let $(\Om\coloneqq\Om(X),\Dd,\Or\coloneqq\Or_X)$ be its oriented de Rham algebra.
The \emph{opposite orientation} on $\CP^2$ induces the opposite orientation $\overline{\Or} = -\Or$ on $\Om$.
The de Rham cohomology algebra $\HDR\coloneqq\HDR(X)$ is isomorphic to the algebra~$V$ for $n=4$ and $\K=\R$ by sending $a\in V$ to the Poincar\'e dual of the canonical projective line $\CP^1\subset \CP^2$.
This isomorphism intertwines the orientations~$\Or_*$ and~$\overline{\Or}_*$ on $\HDR$ with the orientations~$\Or_{1}$ and~$\Or_{-1}$ on~$V$, respectively.
Because~$-1$ is not a square in~$\R$, we conclude that $(\Om,\Dd,\Or_*)$ and $(\Om,\Dd,\overline{\Or}_*)$ are not weakly homotopy equivalent as PDGAs (cf.~Remark~\ref{Rem:Degree}).
\end{example}

\begin{example}[Two small dPD models that are not isomorphic]\label{Ex:SU6}
Consider the Lie group $X\coloneqq \SU(6)$ as an oriented closed manifold of dimension $n\coloneqq 35$, and let $(\Om\coloneqq\Om(X),\Dd,\Or\coloneqq\Or_X)$ be its oriented de Rham algebra. 
The inclusion of the dg-subalgebra of biinvariant forms $(\Om_I,\Dd\equiv 0)\subset(\Om,\Dd)$ is a quasi-isomorphism by \cite[Theorem~1.30]{Felix2008}, so $\HH\coloneqq\Om_I$ is a complement of $\im\Dd$ in $\ker\Dd$ (in $\Om$).
According to \cite[Corollary~1.86]{Felix2008}, the de Rham cohomology algebra $\HDR\coloneqq\HDR(X)$ is isomorphic to the free graded commutative $\R$-algebra $\Lambda_\R(h_3, h_5, h_7, h_9, h_{11})$ generated by cohomology classes $h_i\in \HDR^i$ in degrees $i\in\{3,5,7,9,11\}$.
Pick a representative $\omega_i\in\Omega_I^i$ for each $h_i$, i.e., $[\omega_i]=h_i$.
Given $\xi=(\xi_6,\xi_8)\in \Om^6\times \Om^8$, let $\HH_\xi\coloneqq\oplus_{i=0}^n\HH_{\xi}^i$ be the subspace of $\Om$ obtained from $\HH=\oplus_{i=0}^n \HH^i$ by setting
\[
	\HH_\xi^7 \coloneqq \Span_\R\{\omega_{7,\xi}\coloneqq \omega_7+\Dd\xi_6\},\quad\HH_{\xi}^9 \coloneqq \Span_\R\{\omega_{9,\xi}\coloneqq\omega_9 + \Dd\xi_8\},
\]
and $\HH_\xi^i \coloneqq \HH^i$ for all $i\not\in\{7,9\}$.
Clearly, $\HH_\xi$ is a complement of $\im\Dd$ in $\ker\Dd$.
By Remark~\ref{rem:existencehodge} there are Hodge decompositions $\Om = \HH\oplus \im\Dd \oplus C$ and $\Om=\HH_\xi\oplus\im\Dd\oplus C_\xi$.%
\footnote{Note that a Hodge decomposition $\Om=\Omega_I\oplus\im\Dd\oplus \im\Dd^*$ can also be constructed analytically using Hodge theory for a biinvariant Riemannian metric on $X$.}
Let $\Pp$ and $\Pp_\xi$ be the corresponding special propagators from Lemma~\ref{Lem:HodgeDecompositionsAndPropagators}, and let $\SS_{\Pp}\subset\Om$ and $\SS_{\Pp_\xi}\subset\Om$ be the corresponding small subalgebras from Definition~\ref{Def:SmallSubalgebra}, respectively.
Lemma~\ref{Lem:SmallSubalgebra}\,(a) implies that $\SS_\Pp=\HH$ and that~$\SS_{\Pp_\xi}$ contains the vectors 
\begin{align*}
	\eta_{9,11}&\coloneqq \omega_{9,\xi} \wedge \omega_{11} - \omega_9 \wedge \omega_{11}=\Dd(\xi_8\wedge \omega_{11}),\\
	\eta_{7,9} &\coloneqq \omega_{7,\xi} \wedge \omega_{9,\xi} -\omega_7 \wedge \omega_9 = \Dd(-\omega_7\wedge \xi_8 + \xi_6 \wedge \omega_9 + \xi_6 \wedge \Dd\xi_8),
\end{align*}
and $\Pp_\xi(\eta_{7,9})$.
Using $\Dd \circ \Pp_\xi \circ \Dd = -\Dd$ and Stokes' theorem, we compute the following:
\allowdisplaybreaks
\begin{align*}
	\la \Pp_\xi(\eta_{7,9}),\eta_{9,11}\ra
		&=\int_X \Pp_\xi(\eta_{7,9})\wedge \Dd(\xi_8\wedge\omega_{11}) \\
		&=\int_X (\Dd\circ \Pp_\xi \circ \Dd) (- \omega_7 \wedge \xi_8 + \xi_6\wedge \omega_9  + \xi_6 \wedge \Dd \xi_8)\wedge (\xi_8\wedge \omega_{11})\\
		&=-\int_X (\omega_7 \wedge \Dd \xi_8 + \Dd \xi_6\wedge \omega_9  + \Dd \xi_6 \wedge \Dd \xi_8)\wedge (\xi_8\wedge \omega_{11})\\
		&=-\int_X \xi_8\wedge \Dd \xi_6 \wedge \omega_9 \wedge \omega_{11}.
\end{align*}
This integral can be made nonzero by a suitable choice of $\xi\in\Om^6\times\Om^8$.%
\footnote{Setting $\omega_{9,11}\coloneqq\omega_9\wedge\omega_{11}$, we have $0\neq [\omega_{9,11}]\in\Hh_{\mathrm{dR}}^{20}$, so there exists $p\in X$ such that $\omega_{9,11}(p)\neq 0$.
We can then construct $\xi_6\in\Om^6$ such that $\Dd \xi_6(p)\wedge\omega_{9,11}(p)\neq 0$ in local coordinates around~$p$ and extend it to~$X$ using a partition of unity.
The existence of $\xi_8\in \Om^8$ such that $\la \xi_8, \Dd \xi_6\wedge\omega_{9,11}\ra \neq 0$ now follows from the nondegeneracy of $\la-,-\ra\colon\Om\times\Om\to\R$ (which can be proven by a similar argument).}

The inclusion $\iota\colon S_{\Pp_\xi}\into \Om$ and the quotient map $\pi_\QQ\colon \SS_{\Pp_\xi}\onto \QQ(\SS_{\Pp_\xi})$ are quasi-isomorphisms by Lemmas~\ref{Lem:SmallSubalgebra}\,(c) and~\ref{Lem:HodgeType}\,(a), respectively.
Therefore, $[\omega_{9,11}]\neq 0$ implies $[\pi_\QQ(\omega_{9,11})]\neq 0$, and we have $[\pi_\QQ(\eta_{9,11})]=0$.
Since $\la \Pp_\xi(\eta_{7,9}),\eta_{9,11}\ra\neq 0$, we have $\pi_\QQ(\eta_{9,11})\neq 0$, so the vectors $\pi_\QQ(\omega_{9,11}), \pi_\QQ(\eta_{9,11})\in \QQ^{20}(\SS_{\Pp_\xi})$ must be linearly independent.
We see that $\dim \QQ^{20}(\SS_{\Pp_\xi}) \ge 2$.
Since $\QQ^{20}(\SS_\Pp)=\HH^{20} = \Span_\K\{\omega_{9,11}\}$, the dPD algebras $\QQ(\SS_\Pp)$ and $\QQ(\SS_{\Pp_\xi})$ cannot be isomorphic.
\end{example}

\begin{example}[A weak equivalence of PDGAs that cannot be replaced by a single PDGA quasi-isomorphism]\label{Ex:CP7}	
	Let $X\coloneqq \#^7\CP^2$ be the connected sum of seven copies of $\CP^2$ equipped with the \emph{same} (canonical) orientation.
	Then $X$ is a 1-connected oriented closed $4$-manifold, so the oriented de Rham algebra $(\Om\coloneqq\Om(X),\Dd,\Or\coloneqq\Or_X)$ is an oriented PDGA of degree $n=4$, and the de Rham cohomology $(\HDR\coloneqq\HDR(X),\Or_{*})$ is a 1-connected dPD model of $\Om$ by Remark~\ref{Rem:Formality}.
	A Mayer-Vietoris argument shows that $\HDR\simeq \Span_{\R}\{1, h_1,\dotsc,h_7, h\}$ for some $h_i\in \HDR^2$ and $h\in \HDR^4$ that satisfy $h_i \wedge h_j = \delta_{ij} h$ for all $i,j\in\{1,\dotsc,7\}$.

	Let $(V^\prime,\Dd^\prime,\Or^\prime)$ be a dPD algebra of degree $n=4$.
	Lemma~\ref{Lem:dPDMorphism} implies that any PDGA quasi-isomorphism $\Om\to V^\prime$ or $V^\prime\to\Om$ must be an orientation preserving injection.
	This is clearly impossible for $\Om\to V^\prime$ as $\dim V^\prime<\dim \Om = \infty$.
	As for the opposite direction, we can without loss of generality assume that $V^\prime\to \Om$ is a quasi-isomorphic inclusion $(V^\prime,\Dd^\prime,\Or^\prime)\subset(\Om,\Dd,\Or)$.
	Hence, for each $i\in\{1,\dotsc,7\}$ there exists $\omega_i\in V^{\prime 2}$ such that $[\omega_i]=h_i$, and there exists~$\omega\in V^{\prime 4}$ such that $[\omega]=h$.
	Any non-constant smooth function $f\colon X\to \R$ generates an infinite-dimensional subalgebra $\Span_\R\{ f^k \mid k\in\N_0\}\subset \Om$.
	Therefore, $\dim V^\prime < \infty$ implies $V^{\prime 0} = \Span_{\R}\{1\}$, and Poincar\'e duality gives $V^{\prime 4} = \Span_{\R}\{\omega\}$.
	The relation $h_i\wedge h_j = \delta_{ij} h$ thus implies the relation $\omega_i \wedge \omega_j = \delta_{ij}\omega$ for each $i, j\in\{1,\dotsc,7\}$.
	Because $h\neq 0$, there exists $p\in X$ such that $\omega(p)\neq 0$.
	For all $\lambda_j\in \R$ and $i,j\in\{1,\dotsc,7\}$ we then have 
	\[
		\omega_i(p)\wedge \sum_{j=1}^7 \lambda_j \omega_{j}(p)  = \lambda_{i} \omega(p),
	\]
	so $\omega_1(p), \dotsc, \omega_7(p)\in \Lambda^2 T^*_p X$ must be linearly independent.
	However, this is impossible since
	\[
		\dim(\Lambda^2 T^*_p X)=\binom{4}{2} = 6.
	\]
	We conclude that any zig-zag of PDGA quasi-isomorphisms that connects $(\Om,\Dd,\Or)$ to its dPD model must consist of at least two arrows (cf.~Proposition~\ref{Prop:Existence}).  
\end{example}

\begin{example}[Two 1-connected dPD models that cannot be both mapped into a single 1-connected dPD model by PDGA quasi-isomorphisms]\label{Ex:NonUniqueSimplyConnectedPDModel}

	Let $\Lambda\coloneqq\Lambda_\K(a,b,c)$ be the free graded commutative $\K$-algebra with $\deg a =2$,  $\deg b=3$, and $\deg c = 5$, and let~$\Dd$ be the differential on $\Lambda$ with $\Dd a = \Dd c = 0$ and $\Dd b = a^2$.
	If $\Char(\K)=2$, we additionally assume that $b^2 = c^2 = 0$.
	We have $\Dd(a^k b)=a^{k+2}$ and $\Dd(a^k b c) = a^{k+2} c$ for each $k\in\N_0$, so $\Hh(\Lambda)\simeq\Span_{\K}\{1,a,c,ac\}$.

	Let $(V,\Dd)$ be the quotient of $(\Lambda,\Dd)$ by the dg-ideal generated by $b$ and all elements in degrees $i>n\coloneqq 7$.
	Likewise, let $(V^\prime,\Dd^\prime)$ be the quotient of $(\Lambda,\Dd)$ by the dg-ideal generated by $c-ab$ and all elements in degrees $i>n$.
	We then have the following:
	\begin{align*}
		V &\simeq \Span_\K\{1,a,c,ac\},\\
		V^\prime &\simeq\Span_\K\{1,a,b,a^2=\Dd b,ab=c,a^2 b = ac\}.
	\end{align*}	
	We equip $(V,\Dd)$ and $(V^\prime,\Dd^\prime)$ with the orientations $\Or\colon V\to\K$ and $\Or^\prime \colon V^\prime\to\K$ given by $\Or(ac)=\Or^\prime(ac)=1$, which make them into dPD algebras of degree $n$, respectively.
	The quotient maps $\pi\colon \Lambda \onto V$ and $\pi^\prime\colon \Lambda \onto V^\prime$ are clearly quasi-isomorphisms, and we have $\Or_*^\prime \circ \pi_*^\prime = \Or_*\circ \pi$.
	Therefore, $(V,\Dd,\Or)$ and $(V^\prime,\Dd^\prime,\Or^\prime)$ are weakly homotopy equivalent as PDGAs.

	Suppose that there is a dPD algebra ~$(V^\dprime,\Dd^{\dprime},\Or^{\dprime})$ of degree~$n$ and PDGA quasi-isomorphisms $\iota\colon V \to V^\dprime$ and $\iota^\prime\colon V^\prime\to V^\dprime$.
	Lemma~\ref{Lem:dPDMorphism} implies that $\iota$ and $\iota^\prime$ are injective and preserve chain-level orientation.
	Consider the elements $x \coloneqq \iota(\pi(a))$ and $x^\prime\coloneqq \iota^\prime(\pi^\prime(a))$ in $V^{\dprime 2}$.
	Because we deal with quasi-isomorphism and we have $\Hh^2\simeq\Span_\K\{a\}$, there must exist $\alpha\in \K\backslash\{0\}$ such that $[x^\prime]=\alpha [x]$. 
	However, we have $x^2 = \iota(\pi(a^2)) = 0$ since $\pi(a^2)=0$, whereas we have $x^{\prime 2}=\iota^\prime(\pi^\prime(a^2))\neq 0$ since $\pi^\prime(a^2)=\Dd b\neq 0$.
	Therefore, there must exist a nonzero $y\in V^{\dprime 1}$ such that $\Dd y = x^\prime - \alpha x$, so $V^{\dprime}$ cannot be 1-connected. 
	We conclude that $V$ and $V^\prime$, despite being 1-connected, cannot be both mapped to a single \emph{1-connected} dPD algebra $V^\dprime$ (cf.~Proposition~\ref{Prop:Uniqueness}).
\end{example}

\printbibliography

\end{document}